\documentclass[10pt]{amsart}
\usepackage[dvipdfmx]{graphicx}
\usepackage{amsmath,amscd,amsthm,amsxtra}
\usepackage{epsfig,graphics,color,colortbl}
\usepackage{amssymb,latexsym}
\usepackage{mathrsfs,upgreek}
\usepackage[poly,all]{xy}
\usepackage{hyperref}
\usepackage{tikz-cd}
\usepackage{tabmac}

\renewcommand{\baselinestretch}{1.2}

\newtheorem{thm}{\bf Theorem}[section]
\newtheorem{df}[thm]{\bf Definition}
\newtheorem{prop}[thm]{\bf Proposition}
\newtheorem{cor}[thm]{\bf Corollary}
\newtheorem{lem}[thm]{\bf Lemma}
\newtheorem{rem}[thm]{\bf Remark}
\newtheorem{ex}[thm]{\bf Example}

\numberwithin{equation}{section}

\newcommand{\mc}{\mathcal}
\newcommand{\mf}{\mathfrak}

\newcommand{\pf}{\noindent{\bfseries Proof. }}
\newcommand{\ov}{\overline}

\newcommand{\U}{{\mc U}}
\newcommand{\V}{\mc{V}}
\newcommand{\W}{\mc{W}}
\newcommand{\cP}{\mathscr{P}}
\newcommand{\cO}{\mc{O}}
\newcommand{\I}{\mathbb{I}}
\newcommand{\be}{{\bf e}}

\newcommand{\N}{\mathbb{N}}
\newcommand{\Z}{\mathbb{Z}}
\newcommand{\Q}{\mathbb{Q}}

\newcommand{\e}{\epsilon}
\newcommand{\de}{\delta}
\newcommand{\te}{\widetilde{e}}
\newcommand{\tf}{\widetilde{f}}
\newcommand{\td}{\widetilde}

\newcommand{\La}{\Lambda}

\newcommand{\la}{\lambda}

\newcommand{\nw}{^{\nwarrow}}
\newcommand{\se}{^{\searrow}}

\newcommand{\stirlingii}{\genfrac{\{}{\}}{0pt}{}}
\newcommand{\si}{(-1)^{\e_i}}
\newcommand{\sj}{(-1)^{\e_j}}

\newcommand{\blue}[1]{{\color{blue}#1}}

\newcommand{\ot}{\otimes}

\begin{document}
\title[KR modules of generalized quantum group]
{Kirillov-Reshetikhin modules of generalized quantum group of type $A$}

\author{JAE-HOON KWON}

\address{Department of Mathematical Sciences and Research Institute of Mathematics, Seoul National University, Seoul 08826, Korea}
\email{jaehoonkw@snu.ac.kr}

\author{MASATO OKADO}

\address{Department of Mathematics, Osaka City University, Osaka 558-8585, Japan}
\email{okado@sci.osaka-cu.ac.jp}

\keywords{quantum group, crystal base, Lie superalgebra}
\subjclass[2010]{17B37, 22E46, 05E10}

\thanks{J.-H.K. was supported by the National Research Foundation of Korea(NRF) grant funded by the Korea government(MSIT) (No. 2019R1A2C108483311).
M.O. was supported by Grants-in-Aid for Scientific Research No. 16H03922 and No. 19K03426 from JSPS}

\begin{abstract}
The generalized quantum group of type $A$ is an affine analogue of quantum group associated to a general linear Lie superalgebra, which appears in the study of solutions to the tetrahedron equation or the three-dimensional Yang-Baxter equation. In this paper, we develop the crystal base theory for finite-dimensional representations of generalized quantum group of type $A$. As a main result, we construct Kirillov-Reshetikhin modules, that is, a family of irreducible modules which have crystal bases. We also give an explicit combinatorial description of the crystal structure of Kirillov-Reshetikhin modules, the combinatorial $R$ matrix, and energy function on their tensor products.
\end{abstract}

\maketitle
\setcounter{tocdepth}{1}
\tableofcontents

\section{Introduction}

The crystal base of an integrable highest weight module over the quantum group $U_q(\mf{g})$ for a symmetrizable Kac-Moody algebra $\mf{g}$ introduced by Kashiwara in \cite{Ka91} has been one of the most important objects in studying representations of quantum groups, revealing remarkably rich combinatorial structures of them. 
For an affine quantum group $U_q'(\mf{g})$ without derivation, 
there is also a category of finite-dimensional representations 
where the existence of crystal base is not assured.
However, there is a celebrated family of modules $\W_s^{(r)}$ called Kirillov-Reshetikhin modules
\cite{KR2,CH}, and it is conjectured in \cite{HKOTT02,HKOTY99} that they also have crystal bases, where $r$ is the index of the simple root of the classical subalgebra and $s$ is a positive integer. The existence of crystal base is shown when $\mf{g}$ is of non-exceptional type in \cite{OS08}, and the description of crystal structures is given in \cite{FOS09}. 
The KR crystals form an important class of finite affine crystals having a close connection with the theory of perfect crystals and the path realization of crystals of integrable highest weight modules, and applications to various topics in algebraic combinatorics and mathematical physics including Macdonald polynomials at $t=0$, rigged configuration bijection, combinatorial Bethe Ansatz, box-ball systems and so on (cf.~\cite{KMN2,LNSSS,OSS,KSY,HKOTY02}).

Recently, a new Hopf algebra called the {\em generalized quantum group $\U(\e)$} was introduced in \cite{KOS}. It is associated to an $R$ matrix reduced from $n$-products of operators satisfying the tetrahedron equation. 
The tetrahedron equation is a generalization of the Yang-Baxter equation which plays a key role in quantum integrability in three dimensions (see \cite{KOS} and references therein). 
The parameter $\epsilon=(\e_1,\ldots,\e_n)$ with $\e_i\in\{0,1\}$ for $1\leq i\leq n$ stands for the sequence of two kinds of operators occurring in the $n$-product of operators. Note that $\U(\e)$ is equal to the usual quantum group of affine type when $\e$ is homogeneous, that is, $\epsilon_i=\epsilon_j$ for all $i\neq j$,
while for non-homogeneous $\e$ the quantum group associated to a Lie superalgebra enters in a natural way  interpolating two homogenous cases.

The purpose of this paper is to study the crystal base of finite-dimensional modules over $\U(\e)$ when it is of type $A$. In this case, we may view $\U(\e)$ as an affine analogue of the quantum group associated to general linear Lie superalgebra $\mf{gl}_{M|N}$ (cf.~ \cite{Ya94}), where $M$ and $N$ are the numbers of $0$ and $1$ in $\e$, respectively. We assume that the parameter $\e$ is a standard one $\e_{M|N}=(0^M,1^N)$ with $M,N\geq 1$.
As a main result, we construct a family of finite-dimensional irreducible $\U(\e)$-modules having crystal bases, which we call KR modules (Theorem \ref{thm:crystal base of KR}). We show that there exists an irreducible $\U(\e)$-module $\W_s^{(r)}$ with a crystal base, whose crystal $B^{r,s}$ is realized as the set of $(M|N)$-hook semistandard tableaux of shape $(s^r)$ for each $r, s\geq 1$ such that the rectangular partition $(s^r)$ is an $(M|N)$-hook partition. This is a natural super-analogue of the KR crystals of type $A_{n-1}^{(1)}$ since the crystal of a polynomial $\mf{gl}_{M|N}$-module is parametrized by $(M|N)$-hook semistandard tableaux as shown in \cite{BKK}.

The key ingredients for the construction of $\W_s^{(r)}$ and its crystal base are as follows.
The first one is the $R$ matrix given in \cite{KOS} acting on the two-fold tensor product of the $q$-analogue of  the supersymmetric space $\W_s$ with an explicit spectral decomposition. The representation space $\W^{(r)}_s$ is constructed
by fusion construction in \cite{KMN2} using the $R$ matrix. The next one is a polynomial representation of the subalgebra $\ov\U(\e)$ of finite type, whose crystal coincides with that of a polynomial representation of the quantum group $U_q(\mf{gl}_{M|N})$ in \cite{BKK}. It should be pointed out that we can not simply adopt the results in \cite{BKK} even though we have an isomorphism of algebras between $\ov\U(\e)$ and $U_q(\mf{gl}_{M|N})$ (up to a certain extension), since it does not preserve comultiplication. For this reason we give a self-contained proof of the existence of the crystal base of a polynomial $\ov\U(\e)$-module (Theorem \ref{thm:crystal base of poly repn}).

The crystal structure of $B^{r,s}$ with respect to $\ov\U(\e)$ is a well-known one as given in \cite{BKK} since $\W_s^{(r)}$ is an irreducible $\ov\U(\e)$-module.
On the other hand, it is more non-trivial to describe the Kashiwara operator $\tf_0$ with respect to the additional simple root $\alpha_0$ in $\U(\e)$.
In case of type $A_{n-1}^{(1)}$, the description of $\tf_0$ is given by using so-called the promotion operator based on the cyclic symmetry of the Dynkin diagram of $A_{n-1}^{(1)}$. 
We do not have such symmetry in our case and we use instead the switching algorithm of tableaux in \cite{BSS} to have a combinatorial description of $\tf_0$ on $B^{r,s}$.
We also examine the combinatorial $R$ matrix, that is, an isomorphism of crystals from
$B^{r_1,s_1}\ot B^{r_2,s_2}$ to $B^{r_2,s_2}\ot B^{r_1,s_1}$, 
and the associated energy function. Both of them are important objects in general when we consider one dimensional sum in the study of solvable lattice models (see \cite{O:Memoirs} for more details). We give an explicit description of them in terms of insertion scheme of $(M|N)$-hook semistandard tableaux (cf.~\cite{BR}) (Theorem \ref{th:comb R}). It turns out that their combinatorial rule does not differ from that for $U_q(A^{(1)}_{n-1})$ in \cite{S}.

The paper is organized as follows. In Section 2, we briefly review the definition of generalized quantum group $\U(\e)$ of type $A$ and introduce a category $\mc O_{\geq 0}$ of finite-dimensional $\U(\e)$-modules. In Section 3, we define the notion of crystal base of a module in $\mc O_{\geq 0}$. In Section 4, we prove the existence of the crystal base of a polynomial representation of $\ov\U(\e)$ and discuss a relation with the results in \cite{BKK}. In Section 5, we construct a KR module $\W_s^{(r)}$ by fusion construction and prove the existence of its crystal base. Then we give a combinatorial description of the crystal $B^{r,s}$ of $\W_s^{(r)}$ in Section 6, and study the combinatorial $R$ matrix and the energy function in Section 7.

Throughout the paper, we let $q$ be an indeterminate. We put 
\begin{equation*}
\begin{split}
[m]&=\frac{q^m-q^{-m}}{q-q^{-1}}\quad (m\in \Z_{\geq 0}),\\
[m]!&=[m][m-1]\cdots [1]\quad (m\geq 1),\quad [0]!=1,\\
\begin{bmatrix} m \\ k \end{bmatrix}&= \frac{[m][m-1]\cdots [m-k+1]}{[k]!}\quad (0\leq k\leq m).
\end{split}
\end{equation*}

\vskip 5mm

\noindent {\bf Acknowledgements} The first author would like to thank Department of Mathematics of Osaka City University for its hospitality during his visits in 2017 and 2018. The second author would like to thank Atsuo Kuniba and Sergey Sergeev for the collaboration \cite{KOS} on which this work is partially based.

\section{Generalized quantum group of type $A$}\label{sec:generalized quantum group}

\subsection{Generalized quantum group ${\U}(\e)$}


We fix a positive integer $n\geq 4$ throughout the paper. 
Let $\e=(\e_1,\cdots,\e_n)$ be a sequence with $\e_i\in \{0,1\}$ for $1\leq i\leq n$. We denote by $\I$ the linearly ordered set $\{1<2<\cdots <n\}$ with $\Z_2$-grading given by $\I_0=\{\,i\,|\,\e_i=0\,\}$ and $\I_1=\{\,i\,|\,\e_i=1\,\}$.

Let $P=\bigoplus_{i\in \I}\Z\delta_i$ be the free abelian group generated by $\de_i$ with a symmetric bilinear form $(\,\cdot\,|\,\cdot\,)$ given by $(\de_i|\de_j)=(-1)^{\e_i}\de_{ij}$ for $i,j\in \I$. Let $\{\,\de^\vee_i\,|\,i\in \I\,\}\subset P^\vee:={\rm Hom}_\Z(P,\Z)$ be the dual basis such that $\langle \de_i, \de^\vee_j \rangle =\de_{ij}$ for $i,j\in \I$.

Let $I=\{\,0,\ldots,n-1\,\}$ and
\begin{equation}\label{eq:simple root}
\alpha_i=\de_i-\de_{i+1},\quad \alpha_i^\vee = \de^\vee_i-(-1)^{\e_i+\e_{i+1}}\de^\vee_{i+1} \quad (i\in I). 
\end{equation}
Here we understand the subscript $i$ modulo $n$.
Note that $(\alpha_i|\alpha_i)=\pm 2$ if and only if $\e_i=\e_{i+1}=0,1$ respectively, and 
$(\alpha_i|\alpha_i)=0$ if and only if $\e_i\neq \e_{i+1}$.
Put
\begin{equation*}
I_{\rm even}=\{\,i\in I\,|\,(\alpha_i|\alpha_i)=\pm 2\,\},\quad 
I_{\rm odd}=\{\,i\in I\,|\,(\alpha_i|\alpha_i)=0\,\}.
\end{equation*}
We have 
$\si\langle \la,\alpha_i^\vee \rangle =(\la |\alpha_i)$,
for $\la\in P$ and $i\in I$. In particular, 
$\si\langle \alpha_j,\alpha_i^\vee \rangle =(\alpha_j |\alpha_i)=\sj\langle \alpha_i,\alpha_j^\vee \rangle$ for $i,j\in I$.

Set 
\begin{equation*}
q_i=\si q^{\si}=
\begin{cases}
q & \text{if $\e_i=0$},\\
-q^{-1} & \text{if $\e_i=1$},\\
\end{cases} \quad\quad
q_{ij}=
\begin{cases}
q_j, & \text{if $j=i$},\\
q^{-1}_{j}, & \text{if $j=i+1$},\\
1, & \text{otherwise},
\end{cases}
\end{equation*}
for $i,j\in \I$ with $1\leq i\leq j\le n$.
\begin{df}\label{def:U(e)}
{\rm
We define $\U(\e)$ to be the associative $\Q(q)$-algebra with $1$ 
generated by $q^{h}, e_i, f_i$ for $h\in P^\vee$ and $i\in I$ 
satisfying
{\allowdisplaybreaks
\begin{align}
& q^0=1, \quad q^{h +h'}=q^{h}q^{h'} \hskip 2.5cm  (h, h' \in P^{\vee}),\label{eq:Weyl-rel-1} \\ 
& \omega_je_i\omega_j^{-1}=q_{ij}e_i,\quad \omega_jf_i\omega_j^{-1}=q^{-1}_{ij}f_i, \label{eq:Weyl-rel-2} \\ 
&  e_if_j - f_je_i =\delta_{ij}\frac{\omega_i\omega_{i+1}^{-1} - \omega_i^{-1}\omega_{i+1}}{q-q^{-1}},\label{eq:Weyl-rel-3}\\
& e_i^2= f_i^2 =0 \hskip 4.5cm (i\in I_{\rm odd}),
\end{align}
where $\omega_i=q^{\si\de^\vee_i}$ ($i\in \I$), and the Serre-type relations
\begin{equation}\label{eq:Serre-rel-1}
\begin{split}
&\ \, e_i e_j -  e_j e_i = f_i f_j -  f_j f_i =0,
 \hskip 1cm \text{($|i-j|>1$)},\\ 
&
\begin{array}{ll}
e_i^2 e_j- (-1)^{\e_i}[2] e_i e_j e_i + e_j e_i^2= 0,\\ 
f_i^2 f_j- (-1)^{\e_i}[2] f_i f_j f_i+f_j f_i^2= 0,
\end{array}
\ \hskip 1cm\text{($i\in I_{\rm even}$ and $|i-j|=1$)}, 
\end{split}
\end{equation}
and
\begin{equation}\label{eq:Serre-rel-2}
\begin{array}{ll}
  e_{i}e_{i-1}e_{i}e_{i+1}  
- e_{i}e_{i+1}e_{i}e_{i-1} 
+ e_{i+1}e_{i}e_{i-1}e_{i} \\  
\hskip 2cm - e_{i-1}e_{i}e_{i+1}e_{i} 
+ (-1)^{\e_i}[2]e_{i}e_{i-1}e_{i+1}e_{i} =0, \\ 
  f_{i}f_{i-1}f_{i}f_{i+1}  
- f_{i}f_{i+1}f_{i}f_{i-1} 
+ f_{i+1}f_{i}f_{i-1}f_{i}  \\  
\hskip 2cm - f_{i-1}f_{i}f_{i+1}f_{i} 
+ (-1)^{\e_i}[2]f_{i}f_{i-1}f_{i+1}f_{i} =0,
\end{array}\quad \text{($i\in I_{\rm odd}$)}.
\end{equation}}
Here the subscripts $i,j$ are also understood to be modulo $n$. 
We call $\U(\e)$ the {\em generalized quantum group of affine type $A$ associated to $\e$} (see \cite{KOS,Ma}).
}
\end{df}

Note that
\begin{equation*}
\U(\e) \cong 
\begin{cases}
U_q(A_{n-1}^{(1)}), & \text{if $\e_1=\cdots=\e_n=0$},\\
U_{-q^{-1}}(A_{n-1}^{(1)}), & \text{if $\e_1=\cdots=\e_n=1$},
\end{cases}
\end{equation*}
where $U_q(A_{n-1}^{(1)})$ is the usual quantum affine algebra of type $A_{n-1}^{(1)}$ with $P$ the set of classical weights of level 0.

Put $k_i=\omega_i\omega_{i+1}^{-1}=q^{\si\alpha_i^\vee}$ for $i\in I$ and let
\begin{equation*}
D_{ij}=D_{ji}= \prod_{k\in \{i, i+1\}\cap\{j,j+1\}} q_k^{2\delta_{ij}-1}
\quad (i,j\in I).
\end{equation*}
Then the relations in \eqref{eq:Weyl-rel-2} and \eqref{eq:Weyl-rel-3} imply that
\begin{equation*}
k_ie_jk_i^{-1}=D_{ij}e_j,\quad 
k_if_jk_i^{-1}=D_{ij}^{-1}f_j,\quad
e_if_j - f_je_i =\delta_{ij}\frac{k_i - k_i^{-1}}{q-q^{-1}}.
\end{equation*}
For $i\in I$, we put $\U(\e)_i =\langle e_i, f_i, k_i^{\pm 1} \rangle$. We have
\begin{equation*}
\U(\e)_i\cong U_{q_i}(\mf{sl}_2) \quad (i\in I_{\rm even}).
\end{equation*} 


\begin{ex}{\rm If $n=6$ and $\e=(000111)$, then
\begin{equation*}
(D_{ij})_{0\leq i,j\leq 5}
=
\begin{pmatrix}
\!\!\!-1 & q^{-1} & 1 & 1 & 1 & \!\!\!-q \\
\ \ q^{-1}   &  q^2 & \ \ q^{-1} & 1 & 1 & 1 \\
 1  &  \ \ q^{-1} & q^2 & \ \ q^{-1} & 1 & 1 \\
 1  &  1 & \ \ q^{-1} & \!\!-1 & \!\!\!-q & 1 \\
1  &  1 & 1 &  \!\!\!-q & \ \  \ q^{-2} & \!\!\!-q \\
\!\!\!-q & 1  &  1 & 1 &  \!\!\!-q & \ \  \ q^{-2}
\end{pmatrix}.
\end{equation*}
}
\end{ex}

For $M, N\in\Z_{\geq 0}$ with $M+N=n$, put
$$
\e_{M|N}=(\underbrace{0,\cdots,0}_{M},\underbrace{1,\cdots,1}_{N}).
$$ 
When $\e=\e_{M|N}$, we have $I_{\rm even}=I\setminus \{0, M\}$, $I_{\rm odd}=\{0,M\}$, $\I_0=\{1,\ldots, M\}$ and $\I_1=\{M+1,\ldots,n\}$.
The Dynkin diagram associated to the Cartan matrix $(\langle \alpha_j,\alpha^\vee_i \rangle)_{0\leq i,j\leq n}$ is  
\vskip 4mm
\begin{center}
\hskip -4cm \setlength{\unitlength}{0.25in}
\begin{picture}(24,4)
\put(7.95,1){\makebox(0,0)[c]{$\bigcirc$}}
\put(10.2,1){\makebox(0,0)[c]{$\cdots$}}
\put(12.4,1){\makebox(0,0)[c]{$\bigcirc$}}
\put(16.3,1){\makebox(0,0)[c]{$\bigcirc$}}
\put(18.45,1){\makebox(0,0)[c]{$\cdots$}}
\qbezier(8.1,1.25)(10.4,3.4)(14.2,3.5)
\qbezier(14.8,3.5)(18.45,3.4)(20.45,1.25)
\put(8.2,1){\line(1,0){1.5}}
\put(10.6,1){\line(1,0){1.5}}
\put(12.7,1){\line(1,0){1.5}}
\put(14.8,1){\line(1,0){1.2}}
\put(16.55,1){\line(1,0){1.4}}
\put(18.9,1){\line(1,0){1.4}}
\put(14.5,1){\makebox(0,0)[c]{$\bigotimes$}}
\put(14.5,3.5){\makebox(0,0)[c]{$\bigotimes$}}
\put(20.6,1){\makebox(0,0)[c]{$\bigcirc$}}
\put(14.5,2.8){\makebox(0,0)[c]{\tiny $\alpha_{0}$}}
\put(14.5,0.3){\makebox(0,0)[c]{\tiny $\alpha_{M}$}}
\put(8,0.3){\makebox(0,0)[c]{\tiny $\alpha_{1}$}}
\put(20.6,0.3){\makebox(0,0)[c]{\tiny $\alpha_{n-1}$}}
\put(12.5,0.3){\makebox(0,0)[c]{\tiny $\alpha_{M-1}$}}
\put(16.4,0.3){\makebox(0,0)[c]{\tiny $\alpha_{M+1}$}}
\end{picture}
\end{center}%
\noindent where $\bigotimes$ denotes an isotropic simple root.

There is a Hopf algebra structure on $\U(\e)$, where the comultiplication $\Delta$, the antipode $S$, and the couint $\varepsilon$ are given by
\begin{equation}\label{eq:comult-1}
\begin{split}
\Delta(q^h)&=q^h\otimes q^h, \\ 
\Delta(e_i)&=e_i\otimes 1 + k_i \otimes e_i, \\  
\Delta(f_i)&=f_i\otimes k_i^{-1} + 1 \otimes f_i, \\  
\end{split}
\end{equation}
\begin{gather}
S(q^h)=q^{-h}, \ \ S(e_i)=-e_ik_i^{-1}, \ \  S(f_i)=-k_if_i,\\
\varepsilon(q^h)=1,\ \ \varepsilon(e_i)=\varepsilon(f_i)=0,
\end{gather}
for $h\in P^\vee$ and  $i\in I$.

Let $\eta$ be the anti-automorphism on $\U(\e)$ defined by
\begin{equation}\label{eq:involution eta}
\eta(q^h) = q^h,\quad 
\eta(e_i) = q_i^{-1}f_ik_i,\quad  
\eta(f_i) = q_ik_i^{-1}e_i,
\end{equation}
for $h\in P^\vee$ and $i\in I$. 
It satisfies $\eta^2={\rm id}$ and
\begin{equation}\label{eq:eta preserves comult}
\Delta\circ \eta = (\eta\otimes\eta)\circ \Delta.
\end{equation}

\subsection{Representations of $\U(\e)$}
For a $\U(\e)$-module $V$ and $\mu=\sum_{i}\mu_i\delta_i\in P$, let 
\begin{equation*}
V_\mu 
= \{\,u\in V\,|\,\omega_i u= q_i^{\mu_i} u \ \ (i\in \I) \,\}
= \{\,u\in V\,|\,q^{\de^\vee_i}u= (-1)^{\e_i\mu_i} q^{\mu_i} u \ \ (i\in \I) \,\}
\end{equation*}
be the $\mu$-weight space of $V$. For a non-zero vector $u\in V_\mu$, we denote by ${\rm wt}(u)=\mu$ the weight of $u$. Put ${\rm wt}(V)=\{\,\mu\in P\,|\,V_\mu\neq 0\,\}$.

Let $P_{\geq 0}=\sum_{i\in \I}\Z_{\geq 0}\delta_i$ be the set of {\em polynomial weights}.
Define $\cO_{\geq 0}$ to be the category of $\U(\e)$-modules with objects $V$ such that
\begin{equation}\label{eq:polynomial weight}
V=\bigoplus_{\mu\in P_{\geq 0}}V_\mu \quad \text{with}\  \dim V_\mu < \infty,
\end{equation}
which are closed under taking submodules, quotients, and tensor products.
Note that for $i\in I$ and $\mu=\sum_{i}\mu_i\de_i$, we have
\begin{equation}\label{eq:weight space}
\begin{split}
& k_iu=  q_i^{\mu_i}q_{i+1}^{-\mu_{i+1}}u\quad (u\in V_\mu),\quad
 e_i V_\mu\subset V_{\mu+\alpha_i},\quad f_i V_\mu\subset V_{\mu-\alpha_i}.
\end{split}
\end{equation}
In particular, any irreducible $V\in \cO_{\geq 0}$ is finite-dimensional since ${\rm wt}(V)\subset P_{\geq 0}$.

Let $\mc{W}$ be the $\Q(q)$-space given by
\begin{equation*}
\W=\bigoplus_{m_1,\ldots,m_n}\Q(q)| m_1,\ldots,m_n \rangle, 
\end{equation*}
where $m_i\in \Z_{\geq 0}$ if $\e_i=0$, and $m_i\in \{0,1\}$ if $\e_i=1$. For simplicity, we write $|{\bf m}\rangle = | m_1,\ldots,m_n \rangle$ for ${\bf m}=(m_1,\ldots,m_n)\in \Z_{\geq 0}^n$. We put $|{\bf m}|=m_1+\cdots+m_n$. In particular  $|{\bf 0}\rangle$ for ${\bf 0}=(0,\ldots,0)$ denotes the vacuum vector. For $i\in \I$, put $\be_i=(0,\cdots, 1,\cdots, 0)$ where $1$ appears only in the $i$-th component. 

For $s\in \Z_{\geq 0}$, put 
\begin{equation*}
\W_s = \bigoplus_{|{\bf m}|=s}\Q(q)|{\bf m}\rangle.
\end{equation*}
For a parameter $x$, we denote by $\W_s(x)$ a $\U(\e)$-module $V$, where $V=\W_s$ as a $\Q(q)$-space and the actions of $e_i, f_i, \omega_j$ are given by  
\begin{equation}\label{eq:W_l}
\begin{split}
e_i |{\bf m}\rangle &= x^{\delta_{i,0}}[m_{i+1}]|{\bf m} + \be_{i} -\be_{i+1} \rangle,\\
f_i |{\bf m}\rangle &= x^{-\delta_{i,0}}[m_{i}]|{\bf m} - \be_{i} + \be_{i+1} \rangle,\\
\omega_j |{\bf m}\rangle &= q_j^{m_j} |{\bf m} \rangle,
\end{split}
\end{equation}
for $i\in I$, $j\in \I$, and $|{\bf m}\rangle\in \W_s$ (see \cite[Proposition 3.2]{KOS}). It is clear that $\W_s(x)\in \cO_{\geq 0}$.

\section{Crystal bases}\label{sec:crystal base}

\subsection{Crystal base}
Let us introduce the notion of a crystal base of a $\U(\e)$-module in $\mc{O}_{\geq 0}$, which is based on \cite{BKK}.
\begin{lem}\label{lem:nilpotence}
For any $\epsilon\in \{0,1\}^n$ and $V\in \cO_{\geq 0}$, we have 
\begin{itemize}
\item[(1)] $e_i$ and $f_i$ are locally nilpotent on $V$ for $i\in I_{\rm even}$,

\item[(2)] $e_i V_\mu =f_i V_\mu =0$ for $i\in I_{\rm odd}$ and $\mu\in P_{\geq 0}$ such that $\langle \mu,\alpha_i^\vee \rangle=0$,

\item[(3)] a $\U(\e)_i$-submodule of $V$ is completely reducible for $i\in I_{\rm odd}$.
\end{itemize}
\end{lem}
\pf (1) It is clear since ${\rm wt}(V)\subset P_{\geq 0}$.

(2) We may assume that $V_\mu\neq 0$. If $\mu=\sum_{i}\mu_i\delta_i$, then we have $\mu_i+\mu_{i+1}=\langle \mu,\alpha_i^\vee \rangle =0$. 
Since $\mu\in P_{\geq 0}$, it follows that $\mu_i=\mu_{i+1}=0$.
In particular, $\mu\pm\alpha_i =\mu\pm\de_i\mp\de_{i+1}\not\in P_{\ge 0}$, and so $V_{\mu\pm\alpha_i}=0$.
Since $e_iV_\mu\subset V_{\mu+\alpha_i}$ and $f_iV_\mu\subset V_{\mu-\alpha_i}$,
the assertion follows.


(3) Let $u\in V$ be a non-zero weight vector. It is enough to show that $N=\U(\e)_i u$ is completely reducible. If $\langle{\rm wt}(u),\alpha_i^\vee\rangle=0$, then we have $N=\Q(q)u$ by (2). So we assume that $\langle{\rm wt}(u),\alpha_i^\vee\rangle\neq 0$ with $k_iu=\pm q^a u$ for some $a\in \Z\setminus\{0\}$.\vskip 2mm

{\em Case 1}. Suppose that $e_iu= 0$. Then 
$$(e_if_i-f_ie_i)u = e_if_i u = \frac{k_i-k_i^{-1}}{q-q^{-1}}u=\pm [a]u,$$ 
which implies that $f_i u\neq 0$ and 
$N=\Q(q)u \oplus \Q(q)f_iu$ is irreducible.\vskip 2mm

{\em Case 2}. Suppose that $e_iu\neq 0$. Since 
$e_if_ie_i u= (e_if_i-f_ie_i)e_iu =\pm [a]e_i u\neq 0$, 
we have $e_if_ie_i u\neq 0$ and hence $f_ie_i u\neq 0$. 

If $f_iu=0$, then $f_ie_i u=e_if_iu \mp [a]u=\mp [a]u$. Hence $N=\Q(q)u \oplus \Q(q)e_iu$, which is irreducible.
If $f_iu\neq 0$, then $v:=\pm [a]u + f_ie_i u\neq 0$ and 
$f_iv= \pm[a]f_i u$, $e_if_iu=f_ie_iu\pm[a]u=v$. Hence
$$N=\left(\Q(q)e_iu \oplus \Q(q)f_ie_iu\right)\oplus 
\left( \Q(q) v \oplus \Q(q)f_iv \right)$$
is a sum of two-dimensional irreducible modules.
\qed\vskip 2mm

Let $U_q(\mf{sl}_2)=\langle e,f,k^{\pm 1} \rangle$ be the quantized enveloping algebra of $\mf{sl}_2$. For $m\in\Z_{\geq 0}$, let $V(m)_{\pm 1}$ be the $(m+1)$-dimensional irreducible $U_q(\mf{sl}_2)$-module spanned by $v_a$ for $0\leq a\leq m$, where
\begin{equation*}
\begin{split}
k v_a &=\pm q^{m-2a} v_a,\\
e v_a &=
\begin{cases}
\pm [m-a+1] v_{a-1}, & (1\leq a\leq m), \\
0, & (a=0),
\end{cases}\\
f v_a &=
\begin{cases}
[a+1]v_{a+1}, & (0\leq a\leq m-1), \\
0, & (a=m),
\end{cases}
\end{split} 
\end{equation*}
(cf.~\cite{Jan}).
Suppose that $i\in I_{\rm even}$.
If $(\e_i,\e_{i+1})=(0,0)$, then it is clear that $\U(\e)_i \cong U_q(\mf{sl}_2)$ under 
\begin{equation}\label{eq:sl_2 +}
e\mapsto e_i,\quad f\mapsto f_i,\quad k\mapsto k_i.
\end{equation}
Hence any finite-dimensional $\U(\e)_i$-submodule of $V$ is a direct sum of $V(m)_+$ ($m\in\Z_{\geq 0}$) as a $U_q(\mf{sl}_2)$-module.
On the other hand, if $(\e_i,\e_{i+1})=(1,1)$, then the $\Q(q)$-algebra homomorphism $\xi : U_q(\mf{sl}_2) \longrightarrow \U(\e)_i $ given by
\begin{equation}\label{eq:sl_2 -}
\xi(e) = -e_i,\quad \xi(f) = f_i,\quad \xi(k)=k^{- 1}_i
\end{equation}
is an isomorphism, and hence any finite-dimensional $\U(\e)_i$-submodule of $V$ is a direct sum of $V(m)_{(-1)^m}$ ($m\in\Z_{\geq 0}$) as a $U_q(\mf{sl}_2)$-module.
The induced comultiplication $\Delta^\xi := (\xi^{-1}\otimes \xi^{-1})\circ\Delta\circ \xi$ of $U_q(\mf{sl}_2)$ on a tensor product of polynomial representations is 
\begin{equation}\label{eq:comul sl_2 -}
\begin{split}
&\Delta^\xi(k^{\pm 1})  = k^{\pm 1} \otimes k^{\pm 1},\\
&\Delta^\xi(e)  = e\otimes 1 + k^{-1}\otimes e,\\
&\Delta^\xi(f)  = f\otimes k + 1 \otimes f.
\end{split}
\end{equation}

Now let us define the Kashiwara operators $\te_i$ and $\tf_i$ on $V\in \cO_{\geq 0}$ for $i\in I$.

{\em Case 1}. Suppose that $i\in I_{\rm even}$ and $(\e_i,\e_{i+1})=(0,0)$. We define $\te_i$ and $\tf_i$ on $V$ to be the usual Kashiwara operators on the upper crystal base of $U_q(\mf{sl}_2)$-module $V(m)_+$ ($m\in\Z_{\geq 0}$) under the isomorphism $U_q(\mf{sl}_2)\cong \U(\e)_i$ \eqref{eq:sl_2 +}. 
In other words, if we have $u\in V_\mu$ with $u=\sum_{k\geq 0} f^{(k)}_iu_k$, where $f^{(k)}_i=f_i^k/[k]!$, $l_k=\langle {\rm wt}(u_k),\alpha_i^\vee \rangle$ and $e_iu_k=0$ for $k\geq 0$, then we define
\begin{equation*}
\begin{split}
\te_i u = \sum_{k\geq 1} q^{-l_k+2k-1}f^{(k-1)}_iu_k,\quad 
\tf_i u = \sum_{k\geq 0} q^{l_k-2k-1}f^{(k+1)}_iu_k.
\end{split}
\end{equation*}

{\em Case 2}. Suppose that $i\in I_{\rm even}$ and $(\e_i,\e_{i+1})=(1,1)$. We define $\te_i$ and $\tf_i$ on $V$ to be the usual Kashiwara operators on the lower crystal base of $U_q(\mf{sl}_2)$-module $V(m)_{\pm}$ ($m\in\Z_{\geq 0}$) under the isomorphism $U_q(\mf{sl}_2)\cong \U(\e)_i$ \eqref{eq:sl_2 -}. 
Hence, if we have $u\in V_\mu$ with $u=\sum_{k\geq 0} f^{(k)}_iu_k$, where $e_iu_k=0$ for $k\geq 0$, then  
\begin{equation*}
\begin{split}
\te_i u = \sum_{k\geq 1} f^{(k-1)}_iu_k,\quad 
\tf_i u = \sum_{k\geq 0} f^{(k+1)}_iu_k.
\end{split}
\end{equation*}

{\em Case 3}. Suppose that $i\in I_{\rm odd}$ and $(\e_i,\e_{i+1})=(0,1)$. 
For $u\in V_\mu$, we define
\begin{equation*}
\te_i u =e_i u,\quad \tf_i u =\eta(e_i) u.
\end{equation*}

{\em Case 4}. Suppose that $i\in I_{\rm odd}$ and $(\e_i,\e_{i+1})=(1,0)$. 
For $u\in V_\mu$, we define
\begin{equation*}
\te_i u =\eta(f_i)u,\quad \tf_i u =f_i u.
\end{equation*}

Let $A_0$ be the subring of $\Q(q)$ consisting of $f(q)/g(q)$ with $f(q), g(q)\in \Q[q]$ such that $g(0)\neq 0$.

\begin{df}\label{def:crystal base}
{\rm
Let $V\in \cO_{\geq 0}$ be given. A pair $(L,B)$ is a {\em crystal base of $V$} if it satisfies the following conditions:
\begin{itemize}
\item[(1)] $L$ is an $A_0$-lattice of $V$ and 
$L=\bigoplus_{\mu\in P_{\geq 0}}L_\mu$ where $L_\mu=L\cap V_\mu$,

\item[(2)] $B$ is a signed basis of $L/qL$, that is,  $B={\bf B}\cup (-{\bf B})$ where ${\bf B}$ is a $\Q$-basis of $L/qL$, 

\item[(3)] $B=\bigsqcup_{\mu\in P_{\geq 0}}B_\mu$ where $B_\mu\subset (L /q L)_\mu$

\item[(4)] $\te_i L\subset L$, $\tf_iL\subset L$ and $
\te_i B \subset B\cup\{0\}$, $\tf_i B \subset B\cup\{0\}$ for $i\in I$,

\item[(5)] $\tf_ib=b'$ if and only if $\te_ib'=\pm b$ for $i\in I$ and $b, b'\in B$.
\end{itemize}
}
\end{df}

The $A_0$-lattice $L$ is called a {\em crystal lattice of $V$}.
The set $B/\{\pm 1\}$ is equipped with an $I$-colored oriented graph structure, which we call a {\em crystal (graph) of $V$}, where $b\stackrel{i}{\longrightarrow}b'$ if and only if $\tf_ib=b'$. For $i\in I_{\rm even}$ and $b\in B$, let 
\begin{equation}\label{eq:var_phi_epsilon}
\varepsilon_i(b)=\max\{\,k\,|\,\te_i^k b\neq 0\,\},\quad
\varphi_i(b)=\max\{\,k\,|\,\tf_i^k b\neq 0\,\}.
\end{equation}

Throughout the paper, let us often identify $B$ with the crystal associated with $(L,B)$ when there is no confusion.

For $s\geq 0$, let us set $\W_s=\W_s(1)$ as a $\U(\e)$-module. Put
\begin{equation}\label{eq:crystal base of W(r)}
\begin{split}
\mc{L}_s &=\bigoplus_{|{\bf m}|=s}A_0 | {\bf m} \rangle,\quad
\mc{B}_s =\{\,\pm |{\bf m}\rangle \!\!\! \pmod{q\mc{L}_s}\,|\ |{\bf m}|=s\,\}.
\end{split}
\end{equation}

\begin{prop}  
For $s\in \Z_{\geq 0}$, 
\begin{itemize}
\item[(1)] $(\mc{L}_s,\mc{B}_s)$ is a crystal base of $\W_s$, 

\item[(2)] the crystal $\mc{B}_s$ is connected.
\end{itemize}
\end{prop}
\pf (1) Clearly, Definition \ref{def:crystal base}(1)--(3) hold. 
Let us show that $(\mc{L}_s,\mc{B}_s)$ satisfies Definition \ref{def:crystal base}(4).
The proof of Definition \ref{def:crystal base}(5) is left to the reader.

Let $i\in I$ be given. We understand the subscript $i$ in $\e_i$ modulo $n$.
Given $l\geq 0$ and $0\leq k\leq l$, put $|{\bf m}^{(i)}_{l,k}\rangle = | (l-k)\be_i + k\be_{i+1} \rangle$. 
We may consider only $|{\bf m}^{(i)}_{l,k}\rangle$ since the coefficients of $\be_j$ for $j\not\in\{i,i+1\}$ in any $|{\bf m}\rangle\in \mc{L}_s$ are invariant under $\te_i$ and $\tf_i$.

{\em Case 1}. Suppose that $(\e_i,\e_{i+1})=(0,0)$. 
Given $l\geq 0$, 
we have $\te_i |{\bf m}^{(i)}_{l,0}\rangle=0$ and 
\begin{equation*}
\tf_i^k |{\bf m}^{(i)}_{l,0}\rangle 
= q^{k(l-k)}f_i^{(k)} |{\bf m}^{(i)}_{l,0}\rangle
= q^{k(l-k)}\begin{bmatrix} l \\ k \end{bmatrix} |{\bf m}^{(i)}_{l,k}\rangle
 \in (1+qA_0)|{\bf m}^{(i)}_{l,k}\rangle\quad (0\leq k\leq l).
\end{equation*}
This implies that $\te_i\mc{L}_s\subset \mc{L}_s$ and $\tf_i\mc{L}_s\subset \mc{L}_s$.

{\em Case 2}. Suppose that $(\e_i,\e_{i+1})=(1,1)$. 
Since $m_i, m_{i+1}\in \{0,1\}$ for any $|{\bf m}\rangle$, it is clear that $\te_i\mc{L}_s\subset \mc{L}_s$ and $\tf_i\mc{L}_s\subset \mc{L}_s$.

{\em Case 3}. Suppose that $(\e_i,\e_{i+1})=(0,1)$. 
Consider $|{\bf m}^{(i)}_{l,0}\rangle$ for $l\geq 0$. 
If $l=0$, 
then $\te_i |{\bf m}^{(i)}_{0,0}\rangle =\tf_i |{\bf m}^{(i)}_{0,0}\rangle =0$. 
If $l\geq 1$, then 
$e_i|{\bf m}^{(i)}_{l,0}\rangle =0$ and 
\begin{equation*}
\begin{split}
\tf_i |{\bf m}^{(i)}_{l,0}\rangle 
&=\eta(e_i)|{\bf m}^{(i)}_{l,0}\rangle =q_i^{-1}f_ik_i|{\bf m}^{(i)}_{l,0}\rangle 
=q_i^{-1}q_i^{l}[l]|{\bf m}^{(i)}_{l,1}\rangle\\
&= \left( \frac{1-q^{2l}}{1-q^2} \right) |{\bf m}^{(i)}_{l,1}\rangle
\in (1+qA_0)| {\bf m}^{(i)}_{l,1} \rangle. 
\end{split}
\end{equation*}
Also we have
\begin{equation*}
\te_i\tf_i|{\bf m}^{(i)}_{l,0}\rangle = 
\left( \frac{1-q^{2l}}{1-q^2} \right) |{\bf m}^{(i)}_{l,0}\rangle.
\end{equation*}
Hence $\te_i\mc{L}_s\subset \mc{L}_s$ and $\tf_i\mc{L}_s\subset \mc{L}_s$.

{\em Case 4}. Suppose that $(\e_i,\e_{i+1})=(1,0)$. It is similar to {\em Case 3}.

(2) From the proof of (1), we see that
\begin{equation*}
\mc{B}_s = 
\{\, \pm \tf^{a_1}_{i_1}\cdots\tf^{a_r}_{i_r}|{\bf m}_s\rangle\!\!\!\pmod{q\mc{L}_s} \,|\,
r\geq 0, i_1,\ldots,i_r\in I\setminus\{0\}, a_1,\ldots,a_r \geq 1 \,\}, 
\end{equation*}
where 
\begin{equation*}
|{\bf m}_s\rangle =
\begin{cases}
| s\be_1 \rangle, & \text{if $\e_1=0$},\\
| \be_1+\cdots +\be_s \rangle, & \text{if $\e_1=\cdots=\e_s=1$ and $s\leq n$},\\
| \be_1+\cdots +\be_t + (s-t)\be_{t+1} \rangle, & \text{if $\e_1=\cdots=\e_t=1$ and $\e_{t+1}=0$ with $t< s$}.
\end{cases}
\end{equation*}
In particular, the crystal $\mc{B}_s$ is connected.
\qed

\subsection{Tensor product rule}

Now, we have the following version of tensor product rule of crystal bases for $\U(\e)$-modules in $\mc{O}_{\geq 0}$ with respect to $\Delta$ \eqref{eq:comult-1}, which is slightly different from those  in \cite[Proposition 2.8]{BKK} with respect to $\Delta_\pm$ (see \eqref{eq:comult-2}).  

\begin{prop}\label{prop:tensor product rule}
Let $V_1, V_2\in \cO_{\geq 0}$ be given. Suppose that $(L_i,B_i)$ is a crystal base of $V_i$ for $i=1,2$. 
Then $(L_1\otimes L_2, B_1\otimes B_2)$ is a crystal base of $V_1\otimes V_2$, where $B_1\otimes B_2\subset (L_1/qL_1)\otimes (L_2/qL_2)=L_1\otimes L_2/qL_1\otimes L_2$. 
Moreover, the Kashiwara operators $\te_i$ and $\tf_i$ act on $B_1\otimes B_2$ as follows:
\begin{itemize}
\item[(1)] 
if $i\in I_{\rm even}$ and $(\e_i,\e_{i+1})=(0,0)$, then
{\allowdisplaybreaks
\begin{equation}\label{eq:tensor product rule for even +}
\begin{split}
&\te_i(b_1\otimes b_2)= \begin{cases}
\te_i b_1 \otimes b_2, & \text{if $\varphi_i(b_1)\geq\varepsilon_i(b_2)$}, \\ 
b_1 \otimes \te_i b_2, & \text{if $\varphi_i(b_1)<\varepsilon_i(b_2)$},\\
\end{cases}
\\
&\tf_i(b_1\otimes b_2)=
\begin{cases}
\tf_i b_1 \otimes b_2, & \text{if $\varphi_i(b_1)>\varepsilon_i(b_2)$}, \\
 b_1 \otimes \tf_i b_2, & \text{if $\varphi_i(b_1)\leq\varepsilon_i(b_2)$}, 
\end{cases}
\end{split}
\end{equation}}
\item[(2)] 
if $i\in I_{\rm even}$ and $(\e_i,\e_{i+1})=(1,1)$, then
{\allowdisplaybreaks
\begin{equation}\label{eq:tensor product rule for even -}
\begin{split}
&\te_i(b_1\otimes b_2)= \begin{cases}
 b_1 \otimes \te_ib_2, & \text{if $\varphi_i(b_2)\geq\varepsilon_i(b_1)$}, \\ 
\sigma_i \te_i b_1 \otimes b_2, & \text{if $\varphi_i(b_2)<\varepsilon_i(b_1)$},\\
\end{cases}
\\
&\tf_i(b_1\otimes b_2)=
\begin{cases}
b_1 \otimes \tf_i b_2, & \text{if $\varphi_i(b_2)>\varepsilon_i(b_1)$}, \\
\sigma_i \tf_i b_1 \otimes  b_2, & \text{if $\varphi_i(b_2)\leq\varepsilon_i(b_1)$}, 
\end{cases}
\end{split}
\end{equation}}
where $\sigma_i=(-1)^{({\rm wt}(b_2)|\alpha_i)}$,
\item[(3)] 
if $i\in I_{\rm odd}$ and $(\e_i,\e_{i+1})=(0,1)$, then
{\allowdisplaybreaks
\begin{equation}\label{eq:tensor product rule for odd +}
\begin{split}
\te_i(b_1\otimes b_2)=&
\begin{cases}
\te_i b_1\otimes b_2, & \text{if }\langle {\rm wt}(b_1),\alpha^\vee_i\rangle>0, \\ 
b_1\otimes \te_i b_2, & \text{if }\langle {\rm wt}(b_1),\alpha^\vee_i\rangle=0,
\end{cases}
\\
\tf_i(b_1\otimes b_2)=&
\begin{cases}
\tf_i b_1\otimes b_2, & \text{if }\langle {\rm wt}(b_1),\alpha^\vee_i\rangle >0, \\ 
b_1\otimes \tf_i b_2, & \text{if }\langle {\rm wt}(b_1),\alpha^\vee_i\rangle=0,
\end{cases}
\end{split}
\end{equation}}
\item[(4)] 
if $i\in I_{\rm odd}$ and $(\e_i,\e_{i+1})=(1,0)$, then
{\allowdisplaybreaks
\begin{equation}\label{eq:tensor product rule for odd -}
\begin{split}
\te_i(b_1\otimes b_2)=&
\begin{cases}
\te_i b_1\otimes b_2, & \text{if }\langle {\rm wt}(b_2),\alpha^\vee_i\rangle=0, \\ 
b_1\otimes \te_i b_2, & \text{if }\langle {\rm wt}(b_2),\alpha^\vee_i\rangle>0,
\end{cases}
\\
\tf_i(b_1\otimes b_2)=&
\begin{cases}
\tf_i b_1\otimes b_2, & \text{if }\langle {\rm wt}(b_2),\alpha^\vee_i\rangle =0, \\ 
b_1\otimes \tf_i b_2, & \text{if }\langle {\rm wt}(b_2),\alpha^\vee_i\rangle>0.
\end{cases}
\end{split}
\end{equation}}
\end{itemize}
\end{prop}
\pf Let $L=L_1\otimes L_2$ and $B=B_1\otimes B_2$. It is clear that $(L,B)$ satisfies the conditions (1)-(3) in Definition \ref{def:crystal base}. So it suffices to check the conditions (4) and (5)  for each $i\in I$, and 
the tensor product rules \eqref{eq:tensor product rule for even +}-\eqref{eq:tensor product rule for odd -}.

{\em Case 1}. Suppose that $i\in I_{\rm even}$ and $(\e_i,\e_{i+1})=(0,0)$. 
Then the conditions (4) and (5) in Definition \ref{def:crystal base} and \eqref{eq:tensor product rule for even +} follow from the usual crystal base theory for $U_q(\mf{sl}_2)$ (cf.~\cite{HK}). \vskip 2mm

{\em Case 2}. Suppose that $i\in I_{\rm even}$ and $(\e_i,\e_{i+1})=(1,1)$. 
By \eqref{eq:sl_2 -} and \eqref{eq:comul sl_2 -}, we may apply the tensor product rule for  $V(m)_{\pm 1}$ over $U_q(\mf{sl}_2)=\langle e,f,k^{\pm 1} \rangle$, where the order of tensor product is given in a reverse way.

For $m\in \Z_{\geq 0}$ and $\sigma\in \{\pm 1\}$, let $v_m$ be a highest weight vector in $V(m)_\sigma$. Put
\begin{equation*}
L(m)_\sigma=\bigoplus_{0\leq a\leq m} A_0f^{(a)}v_m,\quad 
B(m)_\sigma=\{\,\pm f^{(a)}v_m \!\!\!\! \pmod{qL(m)_\sigma}\,|\,0\leq a\leq m\,\}.
\end{equation*}
Then $(L(m)_\sigma, B(m)_\sigma)$ is a crystal base of $V(m)_\sigma$ in the sense of Definition \ref{def:crystal base} regarding $U_q(\mf{sl}_2)=\U(\e)_i$ under \eqref{eq:sl_2 -}. 
Let $m_1,m_2\in \Z_{\geq 0}$ and $\sigma_1,\sigma_2\in \{\pm 1\}$ be given.
Following, for example, \cite[Theorem 4.4.3]{HK}, one can check that 
$(L(m_1)_{\sigma_1}\otimes L(m_2)_{\sigma_2},B(m_1)_{\sigma_1}\otimes B(m_2)_{\sigma_2})$ is a crystal base of $V(m_1)_{\sigma_1}\otimes V(m_2)_{\sigma_2}$, and the Kashiwara operators $\td{e}$ and $\td{f}$ with respect to the comultiplication \eqref{eq:comul sl_2 -} act as follows:
\begin{equation*}
\begin{split}
&\te(b_1\otimes b_2)= \begin{cases}
b_1 \otimes \te b_2, & \text{if $\varphi(b_2)\geq\varepsilon(b_1)$}, \\ 
\sigma_2\, \te b_1 \otimes b_2, & \text{if $\varphi(b_2)<\varepsilon(b_1)$},\\
\end{cases}
\\
&\tf(b_1\otimes b_2)=
\begin{cases}
b_1 \otimes \tf b_2, & \text{if $\varphi(b_2)>\varepsilon(b_1)$}, \\
\sigma_2\, \tf  b_1 \otimes  b_2, & \text{if $\varphi(b_2)\leq\varepsilon(b_1)$},
\end{cases}
\end{split}
\end{equation*}
\noindent where $\varepsilon$ and $\varphi$ are defined as in \eqref{eq:var_phi_epsilon}.
This implies the conditions (4) and (5) in Definition \ref{def:crystal base} and 
\eqref{eq:tensor product rule for even -}.\vskip 2mm

{\em Case 3}. Suppose that $i\in I_{\rm odd}$ and $(\e_i,\e_{i+1})=(0,1)$. 
Let $V\in \cO_{\geq 0}$ with a crystal base $(L_0,B_0)$.
For a weight vector $u\in L_0$, we have $u=u_0 + e_iu_1$ by Lemma \ref{lem:nilpotence}, where $f_iu_0=f_iu_1=0$. Note that $k_i u_1 = \pm q^{l} u_1$ for some $l\in \N$, and
\begin{equation}\label{eq:odd Kashiwara}
\begin{split}
\tf_i \te_i u_1 
&= \tf_i e_i u_1 =\eta(e_i)e_i u_1 = q^{-1}f_ik_i e_i u_1 = -q^{-1}k_if_ie_i u_1  \\
&= q^{-1}k_i\frac{k_i-k_i^{-1}}{q-q^{-1}} u_1 = \frac{q^{2l}-1}{q^2-1}u_1\equiv u_1 \pmod{qL_0}.
\end{split}
\end{equation}
Since $\tf_i e_iu = \tf_i e_iu_1 = \frac{q^{2l}-1}{q^2-1}u_1 \in L_0$, we have $u_1\in L_0$ and hence $u_0\in L_0$. 
Furthermore, if $u\equiv b \pmod{qL_0}$ for some $b\in B_0$, then we may assume that $u=e_i^{k}v$ for some $v\in L_0$ and $k=0,1$ such that $f_iv=0$. 

Now, from \eqref{eq:odd Kashiwara} and 
\begin{equation*}
\begin{split}
\Delta(e_i) &= e_i\otimes 1 + k_i\otimes e_i,\\
\Delta(\eta(e_i))
&= \Delta(q^{-1}f_ik_i) = \Delta(q^{-1}f_i)\Delta(k_i)= q^{-1}(f_i\otimes k_i^{-1} + 1\otimes f_i)(k_i\otimes k_i)\\
&= \eta(e_i)\otimes 1 + k_i\otimes \eta(e_i),
\end{split}
\end{equation*}
we see that the pair $(L, B)$ satisfies the conditions (4) and (5) in Definition \ref{def:crystal base} and \eqref{eq:tensor product rule for odd +}.

\vskip 2mm

{\em Case 4}. Suppose that $i\in I_{\rm odd}$ and $(\e_i,\e_{i+1})=(1,0)$. The proof is similar to {\em Case 3.} So we leave it to the reader.
\qed

\subsection{Polarization}

Let $V$ be a $\U(\e)$-module with a symmetric bilinear form $(\ ,\ )$. We call $(\ ,\ )$ a polarization if
\begin{equation}\label{eq:polarization}
(x u , v)=  ( u, \eta(x)v ),
\end{equation}
for $x\in \U(\e)$ and $u,v\in V$.

\begin{df}{\rm
Let $V\in \mc{O}_{\geq 0}$ be given. A crystal base $(L,B)$ of $V$ is called polarizable if it satisfies the following: 
\begin{itemize}
\item[(1)] $V$ has a polarization $(\ ,\ )$ such that $(L,L)\subset A_0$,

\item[(2)] if  $(\ , \ )_0$ is the induced bilinear form on $L/qL$, 
then $(b,b')_0=\delta_{bb'}$ for $b,b'\in B$. In particular, $(\ ,\ )_0$ is positive definite on $L/qL$.
\end{itemize}

}
\end{df}

For $s\in\Z_{\geq 0}$, we define a symmetric bilinear form on $\W_s$ by
\begin{equation}\label{eq:polarization on W(r)}
(|{\bf m}\rangle, |{\bf m'}\rangle) = \delta_{{\bf m}, {\bf m}'}
q^{-\sum_{1\leq i<j\leq n}m_im_j}
\frac{1}{[m_1+\cdots+m_n]!}
\prod_{i=1}^{n}[m_i]!.
\end{equation}

\begin{prop}\label{prop:polarization on W(r)} 
The bilinear form in \eqref{eq:polarization on W(r)} is a polarization, and the crystal base $(\mc{L}_s,\mc{B}_s)$ of $\W_s$ is polarizable with respect to  \eqref{eq:polarization on W(r)}.
\end{prop}
\pf 
Let $i\in I$ be given.
For ${\bf m}=(m_1,\ldots,m_n)\in \Z_{\geq 0}^n$, let ${\bf m}'={\bf m} - \be_i + \be_{i+1}$ where we understand the subscript $i$ modulo $n$. Then we can check that
\begin{equation*}
(|{\bf m'}\rangle, |{\bf m'}\rangle) 
= \frac{[m_{i+1}+1]}{[m_i]}{q^{1-m_i+m_{i+1}}} (|{\bf m}\rangle, |{\bf m}\rangle), 
\end{equation*}
which implies that $(\ ,\ )$ is a polarization.
 
We have 
$(|{\bf m}\rangle, |{\bf m'}\rangle) \in \delta_{{\bf m}, {\bf m}'}(1+qA_0)$
for any $|{\bf m}\rangle, |{\bf m}'\rangle \in \W_s$. 
Hence $(\mc{L}_s,\mc{L}_s)\subset A_0$ and the induced bilinear form $(\ , \ )_0$ is positive definite on $\mc{L}_s/q\mc{L}_s$. 
Therefore, $(\mc{L}_s,\mc{B}_s)$ is polarizable with respect to \eqref{eq:polarization on W(r)}.
\qed

\begin{lem}
Let $V_1, V_2$ be $\U(\e)$-modules with polarizations $(\ , \ )_1$ and $(\ , \ )_2$, respectively. Then $V_1\otimes V_2$ has a polarization given by
\begin{equation*}
(u_1\otimes u_2,v_1\otimes v_2) = (u_1,v_1)_1(u_2,v_2)_2
\end{equation*}
for $u_1,v_1\in V_1$ and $u_2,v_2\in V_2$.
\end{lem}
\pf It follows from \eqref{eq:eta preserves comult}. \qed

\begin{prop}\label{prop:polarization}
We have the following.
\begin{itemize}
\item[(1)] If $(L_i,B_i)$ is a polarizable crystal base of $V_i\in \mc{O}_{\geq 0}$ for $i=1,2$, then $(L_1\otimes L_2, B_1\otimes B_2)$ is a polarizable crystal base of $V_1\otimes V_2$.

\item[(2)] If $V\in \mc{O}_{\geq 0}$ has a polarizable crystal base, then $V$ is completely reducible.
\end{itemize}
\end{prop}
\pf The proofs of (1) and (2) are the same as in \cite[Lemma 2.11]{BKK} and \cite[Theorem 2.12]{BKK}, respectively.
\qed

\section{Polynomial representations of finite type $A$}\label{sec:polynomial representation}
\label{sec:polynomial repn}
\subsection{Polynomial representations}

Let $\ov{I}=I\setminus\{0\}$ and let $\ov\U(\e)$ be the subalgebra of ${\U}(\e)$ generated by $q^h$, $e_i, f_i$ for $h\in P^\vee$ and $i\in \ov I$. 
Put $\ov I_{\rm even}=\ov I\cap I_{\rm even}$ and $\ov I_{\rm odd}=\ov I\cap I_{\rm odd}$.

We define $\ov\cO_{\geq 0}$ to be the category of $\ov\U(\e)$-modules $V$ satisfying \eqref{eq:polynomial weight}. We call $V\in\ov\cO_{\geq 0}$ a {\em polynomial representation of $\ov\U(\e)$}. The category $\ov\cO_{\geq 0}$ is closed under taking submodules, quotients, and tensor products.
We define a crystal base of $V\in \ov\cO_{\geq 0}$ in the same way as in Definition \ref{def:crystal base}.

Let $\V=\W_1$ as a $\ov{\U}(\e)$-module with $\Q(q)$-basis 
$\left\{\,|\be_i\rangle  \,|\, i\in \I\,\right\}$.
Clearly $\V\in \ov\cO_{\geq 0}$, and 
$(\mc{L},\mc{B}):=(\mc L_1,\mc B_1)$ is a polarizable crystal base by Proposition \ref{prop:polarization on W(r)}. 
Furthermore, for $m\geq 1$, $(\mc{L}^{\otimes m},\mc{B}^{\otimes m})$ is a polarizable crystal base of $\mc{V}^{\otimes m}$ by Proposition \ref{prop:polarization}.

\begin{cor}\label{cor:complete reducibility of tensor powers}
For $m\geq 1$, $\mc{V}^{\otimes m}$ is completely reducible. 
\end{cor}

\begin{rem}{\rm
Note that
\begin{equation}\label{eq:two extremes}
\ov\U(\e)\cong 
\begin{cases}
U_{q}(\mf{gl}_n), & \text{if $\e_1=\cdots =\e_n=0$},\\ 
U_{-q^{-1}}(\mf{gl}_n), & \text{if $\e_1=\cdots =\e_n=1$},
\end{cases}
\end{equation}
where $U_q(\mf{gl}_n)$ is the quantized enveloping algebra of the general linear Lie algebra $\mf{gl}_n$. 
One may think of the usual crystal base of a $\ov\U(1,\dots ,1)$-module $V$ with respect to $v=-q^{-1}$ and hence a crystal at $v=0$, while the crystal base of $V$ in Definition \ref{def:crystal base} is the one at $q=0$. In this case, it is not difficult to see that each irreducible in $\ov\cO_{\geq 0}$ has a crystal base in the sense of Definition \ref{def:crystal base}, whose crystal is isomorphic to that of usual crystal with the same highest weight.}
\end{rem}

\subsection{Quantum superalgebra $U_q(\mf{gl}_{M|N})$}\label{subsec:Yamane's quantum group}

Let $\Theta$ be the bialgebra over $\Q(q)$ generated by $\vartheta_j$ for $j\in \I_1$, which commute and satisfy $\vartheta_j^2=1$.
Here the comultiplication is given by $\Delta(\vartheta_j)=\vartheta_j\otimes \vartheta_j$ for $j\in \I_1$.
Then $\Theta$ acts on $\ov\U(\e)$ by 
\begin{equation}\label{eq:vartheta-rel}
\begin{split}
&\vartheta_j q^h =q^h,\quad
\vartheta_je_i=(-1)^{(\delta_j|\alpha_i)}e_i,\quad 
\vartheta_jf_i=(-1)^{(\delta_j|\alpha_i)}f_i,
\end{split}
\end{equation}
for $j\in \I_1$, $h\in P^\vee$, and $i\in \ov I$. It is clear that $\ov \U(\e)$ is a $\Theta$-module algebra.
Let $\ov\U(\e)[\vartheta]$ be the semidirect product of $\ov \U(\e)$ and $\Theta$.
Then it is straightforward to check the following.


\begin{prop}\label{prop:action of vartheta}
For $V\in \ov\cO_{\geq 0}$, $V$ can be extended to a $\ov\U(\e)[\theta]$-module, where 
$\vartheta_j u = (-1)^{\mu_j} u$ 
for $j\in \I_1$ and $u\in V_\mu$ with $\mu=\sum_i\mu_i\delta_i$.
\end{prop} 
\qed

Let $U(\e_{M|N})$ be the associative $\Q(q)$-algebra with $1$ 
generated by $q^{h}, E_i, F_i$ for $h\in P^\vee$ and $i\in \ov I$ 
satisfying
{\allowdisplaybreaks
\begin{align*}
& q^0=1, \quad q^{h +h'}=q^{h}q^{h'} \hskip 2.5cm  (h, h' \in P^{\vee}), \\ 
& \varpi_jE_i\varpi_j^{-1}=q^{(\alpha_i|\delta_j)}E_i,\quad 
\varpi_jF_i\varpi_j^{-1}=q^{-(\alpha_i|\delta_j)}F_i, \\ 
&  E_iF_j - (-1)^{p(i)p(j)} F_jE_i =
{(-1)^{\e_i}}\delta_{ij}\frac{\varpi_i\varpi_{i+1}^{-1} - \varpi_i^{-1}\varpi_{i+1}}{q-q^{-1}},\\
& E_M^2 = F_M^2=0,\\
& E_i E_j - (-1)^{p(i)p(j)} E_j E_i = F_i F_j - (-1)^{p(i)p(j)}  F_j F_i =0,
 \hskip 1cm \text{($|i-j|>1$)},\\ 
&\!\!\!
\begin{array}{ll}
E_i^2 E_j- [2] E_i E_j E_i + E_j E_i^2= 0,\\ F_i^2 F_j-[2] F_i F_j F_i+F_j F_i^2= 0,
\end{array}
\ \hskip 1cm\text{($i\not =M$ and $|i-j|=1$)}, \\
&\!\!\!
\begin{array}{ll}
  E_{M}E_{M-1}E_{M}E_{M+1}  
+ E_{M}E_{M+1}E_{M}E_{M-1} 
+ E_{M+1}E_{M}E_{M-1}E_{M} \\  
\hskip 3cm + E_{M-1}E_{M}E_{M+1}E_{M} 
- [2]E_{M}E_{M-1}E_{M+1}E_{M} =0, \\   
  F_{M}F_{M-1}F_{M}F_{M+1}  
+ F_{M}F_{M+1}F_{M}F_{M-1} 
+ F_{M+1}F_{M}F_{M-1}F_{M}  \\  
\hskip 3cm + F_{M-1}F_{M}F_{M+1}F_{M} 
- [2]F_{M}F_{M-1}F_{M+1}F_{M} =0, 
\end{array}
\end{align*}
\noindent where $p(i)=\e_i+\e_{i+1}$ and
$\varpi_j=q^{\sj\de^\vee_j}$ for $j\in \I$. Here the subscripts $i,j$ are also understood to be modulo $n$.}

We also define $U(\e_{M|N})[\vartheta]$ to be the semidirect product of $\ov \U(\e)$ and $\Theta$, where $\vartheta_j$ acts on $U(\e_{M|N})$ as in \eqref{eq:vartheta-rel} with $e_i$ and $f_i$ replaced by $E_i$ and $F_i$ for $i\in \ov I$.  
 
\begin{prop}\label{prop:iso phi}
For $M, N\in\Z_{\geq 0}$ with $M+N=n$, there exists an isomorphism of $\Q(q)$-algebras
$$\tau : U(\e_{M|N})[\vartheta] \longrightarrow \ov\U(\e_{M|N})[\vartheta],$$ 
such that 
{\allowdisplaybreaks
\begin{equation}\label{eq:iso phi}
\begin{split}
& \tau(\vartheta_j) =\vartheta_j \quad (M+1\leq j\leq n),\\ 
& \tau(\varpi_j)= 
\begin{cases}
\omega_j, & \text{$(1\leq j\leq M)$},\\ 
\omega_j\vartheta_j, & \text{$(M+1\leq j\leq n)$},
\end{cases}\\
& \tau(E_i)= 
\begin{cases}
e_i, & \text{$(1\leq i\leq M)$},\\
-e_i(\vartheta_i\vartheta_{i+1})^{i+M}, & \text{$(M+1\leq i\leq n-1)$},
\end{cases}\\
& \tau(F_i)= 
\begin{cases}
f_i, & \text{$(1\leq i< M)$},\\
f_M\vartheta_{M+1}, & \text{$(i=M)$},\\
f_i(\vartheta_i\vartheta_{i+1})^{i+M+1}, & \text{$(M+1\leq i\leq n-1)$}.
\end{cases}
\end{split}
\end{equation}}
\end{prop} 
\pf It is straighforward to check that $\tau(q^h)$, $\tau(E_i)$, $\tau(F_i)$, $\vartheta_j$ satisfy the defining relations of $U(\e_{M|N})[\vartheta]$. 
For example, we have 
{\allowdisplaybreaks
\begin{align*}
&\tau(E_M)\tau(F_M) + \tau(F_M)\tau(E_M) \\
&= e_M f_M \vartheta_{M+1} +  f_M \vartheta_{M+1} e_M = e_M f_M \vartheta_{M+1} -  f_M e_M \vartheta_{M+1} \\ 
&= (e_M f_M  -  f_M e_M) \vartheta_{M+1} 
= \frac{k_M-k_{M}^{-1}}{q-q^{-1}} \vartheta_{M+1}\\ 
&=\frac{{\omega_M}\,{\omega_{M+1}}^{-1} - {\omega_M}^{-1}{\omega_{M+1}}}{q-q^{-1}}\vartheta_{M+1} = \frac{\tau(\varpi_M)\,\tau(\varpi_{M+1})^{-1} - \tau({\varpi_M})^{-1}\tau({\varpi_{M+1}})}{q-q^{-1}},
\end{align*}}
and  
{\allowdisplaybreaks
\begin{align*}
&\tau({E}_i)\tau({F}_i) - \tau({F}_i)\tau({E}_i) \\
&= 
-e_i (\vartheta_i\vartheta_{i+1})^{i+M} f_i(\vartheta_i\vartheta_{i+1})^{i+M+1}
+ f_i(\vartheta_i\vartheta_{i+1})^{i+M+1} e_i (\vartheta_i\vartheta_{i+1})^{i+M} \\ 
&= -(e_i f_i  -  f_i e_i) \vartheta_i\vartheta_{i+1} = -\frac{k_i-k_{i}^{-1}}{q-q^{-1}} \vartheta_i\vartheta_{i+1} \\
&= -\frac{\omega_i\omega_{i+1}^{-1} - \omega_i^{-1}\omega_{i+1}}{q-q^{-1}}\vartheta_i\vartheta_{i+1} = - \frac{\tau({\varpi_i})\tau({\varpi_{i+1}})^{-1} - \tau({\varpi_i})^{-1}\tau({\varpi_{i+1}})}{q-q^{-1}},
\end{align*}}
for $M+1\leq i\leq n-1$. We leave the verification of the other relations to the reader.

So $\tau$ is a well-defined surjective $\Q(q)$-algebra homomorphism. 
Conversely, there is a $\Q(q)$-algebra homomorphism $\tau' : \ov\U(\e_{M|N})[\vartheta] \longrightarrow U(\e_{M|N})[\vartheta]$ similar to $\tau$, which is defined by replacing $E_i$ and $F_i$ with $e_i$ and $f_i$ for $i\in \ov I$, respectively in \eqref{eq:iso phi}. Then $\tau$ and $\tau'$ are inverses of each other. Hence $\tau$ is an isomorphism. 
\qed\vskip 2mm

The operator $\theta_i:=\vartheta_i\vartheta_{i+1}$ for $M\leq i \leq n$ is equal to $\theta_i$ in \cite[Section 3.3.1]{KOS}, where we assume that $\vartheta_M=1$.
If we put $K_i=\varpi_i\varpi_{i+1}^{-1}\in U(\e_{M|N})$ for $i\in \ov I$, then  
\begin{equation}\label{eq:sign twist}
\tau(K_i)=
\begin{cases}
k_i, & (1\leq i<M),\\
k_i\theta_{i}, & (M\leq i\leq n-1).
\end{cases}
\end{equation}
The isomorphism $\tau$ when restricted to the subalgebra generated by $E_i$, $F_i$, $K_i^{\pm 1}$ and $\theta_j$ for $i\in \ov I$ and $M\leq j\leq n$ is equal to the one in \cite[Section 3.3]{KOS}.  
 
\begin{rem}{\rm
For arbitrary $\epsilon$, a quantum superalgebra of type $A$ is defined in \cite{Ya94}. 
Similar to Proposition \ref{prop:iso phi}, this algebra was shown to be isomorphic to $\ov{\U}(\e)$ by adding generators of order 4 in \cite{Ma}.
}
\end{rem}

Set $$\chi=\vartheta_{M}\cdots\vartheta_n.$$ Then $\chi^2=1$, $\chi q^h =q^h \chi$, and  
$\chi X_i = (-1)^{p(i)} X_i\chi$ for $i\in \ov I$ and $X=E, F$. 
We define
\begin{equation*}
U_q(\mf{gl}_{M|N}) = U(\e_{M|N}) \oplus U(\e_{M|N})\chi,
\end{equation*}
which is the subalgebra of $U(\e_{M|N})[\vartheta]$ generated by $U(\e_{M|N})$ and $\chi$. It is isomorphic to the quantum superalgebra associated to the general linear Lie superalgebra $\frak{gl}_{M|N}$ in \cite{Ya94} (see also \cite{BKK}).  

Let $O_{\rm int}$ be the category of $U_q(\mf{gl}_{M|N})$-modules satisfying the following conditions \cite[Definiton 2.2]{BKK}:
\begin{itemize}
\item[(1)] $V=\bigoplus_{\mu\in P}V_\mu$, where 
$V_\mu = \{\,u\in V\,|\,q^{h}u= q^{\langle \mu,h\rangle} u \ \ (h\in P^\vee) \,\}$,

\item[(2)] $\dim V_\mu < \infty$ for $\mu\in P$,

\item[(3)] $V$ is a locally finite $U_q(\mf{gl}_{M|N})_i$-module for $i\neq M$, 
where $U_q(\mf{gl}_{M|N})_i=\langle E_i, F_i, K_i^{\pm 1} \rangle$,

\item[(4)] $V_\mu\neq 0$ implies $\langle \mu, \alpha^\vee_M\rangle\geq 0$,

\item[(5)] $E_M V_\mu =F_M V_\mu =0$ if $\langle \mu,\alpha_M^\vee \rangle=0$ for $\mu\in P$.

\end{itemize}

For $V\in \ov\cO_{\geq 0}$,
let $V^\tau$ be the pullback of $V$ via $\tau$, which is the same as $V$ as a space, whose elements are denoted by $u^\tau$ for $u\in V$, and the action of $U_q(\mf{gl}_{M|N})$ is given by $x u^\tau = (\tau(x)u)^\tau$ for $x\in U_q(\mf{gl}_{M|N})$ and $u\in M$. 


\begin{prop}\label{prop:O_int in BKK}
For $V\in \ov\cO_{\geq 0}$, we have $V^\tau \in O_{\rm int}$.
\end{prop}
\pf 
Property (1) follows from Proposition \ref{prop:action of vartheta} and \eqref{eq:sign twist}.
Properties (2) and (4) are clear.
Properties (3) and (5) follow from Lemma \ref{lem:nilpotence}.
\qed
\vskip 2mm

\begin{rem}\label{rem:comultiplication}
{\rm
There are two Hopf algebra structures on $U_q(\mf{gl}_{M|N})$ with the respective comultiplications $\Delta_\pm$ given by
\begin{equation}\label{eq:comult-2}
\begin{split}
\Delta_\pm(\chi)=\chi\otimes \chi, \quad
\Delta_\pm(q^h) =q^h\otimes q^h, \\ 
\begin{cases}
\Delta_+(E_i) =E_i\otimes 1+ \chi^{p(i)} K_i\otimes E_i,\\ 
\Delta_-(E_i) =E_i\otimes K_i^{-1} + \chi^{p(i)}\otimes E_i, \\  
\end{cases}\\
\begin{cases}
\Delta_+(F_i) =F_i\otimes K_i^{-1} + \chi^{p(i)}\otimes F_i, \\  
\Delta_-(F_i) =F_i\otimes 1+ \chi^{p(i)} K_i\otimes F_i, 
\end{cases}
\end{split}
\end{equation}
for $h\in P$ and $i\in \ov I$.
On the other hand, the comultiplication $\Delta^\tau=:(\tau^{-1}\otimes \tau^{-1})\circ\Delta\circ\tau$ on $U_q(\mf{gl}_{M|N})$ induced from $\Delta$ via $\tau$ is given by
\begin{equation}\label{eq:comult-3}
\begin{split}
&\Delta^\tau(\chi)=\chi\otimes \chi, \ \ 
\Delta^\tau(q^h) =q^h\otimes q^h, \\ 
&\Delta^\tau(E_i) =
\begin{cases}
E_i\otimes 1+ K_i\otimes E_i, & (1\leq i\leq M),\\ 
E_i\otimes \theta_i^{M+i}+  K_i\theta_i^{M+i-1}\otimes E_i, & (M+1\leq i\leq n-1),
\end{cases}\\
&\Delta^\tau(F_i) =
\begin{cases}
F_i\otimes K_i^{-1} + 1 \otimes F_i, & (1\leq i\leq M),\\ 
F_i\otimes K_i^{-1}\theta_i^{M+i-1}+  \theta_i^{M+i}\otimes F_i, & (M+1\leq i\leq n-1), \\  
\end{cases}
\end{split}
\end{equation}
for $h\in P$ and $i\in \ov I$. 

In the next subsection, we discuss the crystal bases of $\ov\U(\e_{M|N})$-modules $V$ in $\ov\cO_{\geq 0}$. By Proposition \ref{prop:O_int in BKK}, one may take a crystal base of $V$ as the one of $U_q(\mf{gl}_{M|N})$-module $V^\tau$ given in \cite{BKK}. But the construction of crystal bases in \cite{BKK} and their tensor products are given with respect to the comultiplication $\Delta_-$ \eqref{eq:comult-2}, and we have $\Delta^\tau \neq \Delta_-$ from \eqref{eq:comult-2} and \eqref{eq:comult-3}. So we give instead a self-contained exposition on the crystal bases of $\ov\U(\e_{M|N})$-modules $V$ in $\ov\cO_{\geq 0}$, though the proofs are essentially the same as those in \cite{BKK}.

}
\end{rem}

\subsection{Crystal bases of polynomial representations}\label{subsec:crystal of poly repn}

Let $\ov\U(\e)^+$ be the subalgebra of $\ov\U(\e)$ generated by $q^h$ and $e_i$ for $h\in P^\vee$ and $i\in \ov I$. For $\la=\sum\la_i\delta_i\in P_{\geq 0}$, let ${\bf 1}_\la=\Q(q)v_\la$ be the one-dimensional $\ov\U(\e)^+$-module such that 
$e_i v_\la =0$ and $\omega_j v_\la = q_j^{\la_j}v_\la$ for $i\in \ov I$ and $j\in \I$. We define $V(\la)$ to be the irreducible quotient of $\ov\U(\e)\otimes_{\ov\U(\e)^+} {\bf 1}_\la$. By \eqref{eq:weight space}, we see that $V(\la)=\bigoplus_{\mu}V(\la)_\mu$, where the sum is over $\mu = \la-\sum_{i\in \ov I}c_i\alpha_i$ with $c_i\in\Z_{\geq 0}$ and ${\rm dim}V(\la)_\mu<\infty$.

Let $\cP$ be the set of all partitions. 
For $M, N\in\Z_{\geq 0}$ with $M+N=n$, a partition $\la=(\la_i)_{i\ge 1}\in \cP$ is called an $(M|N)$-hook partition if $\la_{M+1}\leq N$ (cf.~\cite{BR}). We denote the set of all $(M|N)$-hook partitions by $\cP_{M|N}$.
For example, \vskip 2mm
$$\la=(5,3,2,2)=
\resizebox{.25\hsize}{!}
{\def\lr#1{\multicolumn{1}{|@{\hspace{.75ex}}c@{\hspace{.75ex}}|}{\raisebox{-.04ex}{$#1$}}}
\def\l#1{\multicolumn{1}{|@{\hspace{.75ex}}c@{\hspace{.75ex}}}{\raisebox{-.04ex}{$#1$}}}\raisebox{-.6ex}
{$\begin{array}{cccccccc}
\cline{1-8}
\lr{\ \,} & \lr{\ \, } & \lr{\ \,} & \lr{\ \,} & \lr{\ \,} & & &  \\
\cline{1-5}
\lr{\ }  & \lr{\ } & \lr{\ } &  &  \\
\cline{1-3}
\lr{\ }  & \lr{\ } & & & \\
\cline{1-2}\cline{5-8}
\lr{\ }  & \lr{\ } &  & & \l{\ } \\
\cline{1-2}
\l{\ }  &   &  & &  \l{\ } \\
\end{array}$}}\quad \in \cP_{3|4}
$$

For $\la\in \cP_{M|N}$, we define
\begin{equation}\label{eq:highest weight correspondence}
\La_\la = \la_1\de_1 +\cdots \la_M\de_M + \mu_1\de_{M+1}+\cdots+\mu_N\de_{M+N},
\end{equation}
where $(\mu_1,\ldots,\mu_N)$ is the conjugate of the partition $(\la_{M+1},\la_{M+2},\ldots)$. Since the map $\la \mapsto \La_\la$ is injective, we identify $\cP_{M|N}$ with its image in $P_{\geq 0}$ if there is no confusion.

\begin{prop}\label{prop:irr poly rep}
Let $V$ be an irreducible $\ov\U(\e_{M|N})$-module in $\ov\cO_{\geq 0}$. Then $V\cong V(\la)$ for some $\la\in \cP_{M|N}$. Moreover, $V^\tau$ is an irreducible $U_q(\mf{gl}_{M|N})$-module in $O_{\rm int}$ with highest weight $\la$.
\end{prop}
\pf Since $V$ is finite-dimensional, there exists a highest weight vector $v$ with weight $\la\in P_{\geq 0}$. Let us consider $V^{\tau}$ (see Section \ref{subsec:Yamane's quantum group}).
By Proposition \ref{prop:iso phi}, we see that the set of weights in $V^\tau$ is equal to that of $V$.
Suppose that $V^\tau$ is not irreducible. Then there exists a vector $w$ with weight $\mu\in  \la-\sum_{i\in \ov I}\Z_{\geq 0}\alpha_i$ such that $\mu\neq \la$ and $E_i w=0$ for $i\in \ov I$. This implies that $w$ also generates a proper $\ov\U(\e_{M|N})$-submodule of $V$, which is a contradiction. So $V^\tau$ is an irreducible $U_q(\mf{gl}_{M|N})$-module with highest weight $\la$.
By \cite[Propositions 3.4 and 4.5]{BKK}, it follows that $\la\in\cP_{M|N}$. Hence $V\cong V(\la)$.
\qed


 

\vskip 2mm

Suppose that $\e$ is arbitrary.
Let $\mc{B}$ be the crystal of the natural representation $\mc{V}$ of $\ov{\U}(\e)$.
For simplicity, we identify $\mc{B}$ with $\I$ as a set, 
and $\mc{B}^{\otimes m}$ with a word of length $m$ with letters in $\I$ for $m\geq 1$. 
The crystal $\mc{B}$ is given by
\begin{equation*}
1\ \stackrel{^{1}}{\longrightarrow}\ 2\ \stackrel{^{2}}{\longrightarrow}
\cdots\stackrel{^{n-1}}{\longrightarrow}\  n.
\end{equation*}
Thanks to Proposition \ref{prop:tensor product rule}, the crystal structure on $\mc{B}^{\otimes m}$ coincides with the one in \cite[Section 2]{BKK}.
For the reader's convenience, let us briefly describe $\tf_i$ ($i\in \ov I$) on $\mc{B}^{\otimes m}$.
Let $b=b_1\cdots b_m\in \mc{B}^{\otimes m}$ be given, where $b_t\in \I$ for $1\leq t\leq m$.

\begin{itemize}
\item[{\em Case 1}.] Suppose that $i\in \ov I_{\rm even}$ with $(\e_i,\e_{i+1})=(0,0)$. 
For $1\leq t\leq m$, we put $a_t=+$, $-$, and $\,\cdot\,$ if $b_t=i$, $i+1$, and otherwise, respectively. In the sequence ${\bf a}=(a_1,\cdots,a_m)$, we replace a
pair $(a_{t},a_{t'})=(+,-)$, where $t<t'$ and
$a_{t''}=\cdot$ for $t<t''<t'$, with $(\,\cdot\,,\,\cdot\,)$, and repeat
this process as far as possible until we get a sequence with no $-$
placed to the right of $+$. 
We denote  this sequence by $\td{\bf a}$.
If $t$ corresponds to the leftmost $+$ in $\td{\bf a}$,
then $\tf_i b$ is given by replacing $b_t=i$ with $i+1$. If there is no $+$ in $\td{\bf a}$, then $\tf_ib=0$.

\item[{\em Case 2}.] Suppose that $i\in \ov I_{\rm even}$ with $(\e_i,\e_{i+1})=(1,1)$. In this case, $\tf_i b$ is obtained by the same way as in {\em Case 1}, except $\bf a$ is replaced by ${\bf a}^{\rm rev}:=(a_m,\cdots,a_1)$.

\item[{\em Case 3}.] Suppose that $i\in \ov I_{\rm odd}$ with $(\e_i,\e_{i+1})=(0,1)$. Choose the smallest $1\leq t\leq m$ such that $({\rm wt}(\be_{b_t})|\alpha_i)\neq 0$ or $b_t\in \{i, i+1\}$. If $b_t=i$, then $\tf_i b$ is given by replacing $b_t=i$ with $i+1$. If $b_t=i+1$ or there is no such $t$, then $\tf_ib=0$.

\item[{\em Case 4}.] Suppose that $i\in \ov I_{\rm odd}$ with $(\e_i,\e_{i+1})=(1,0)$. Then $\tf_i b$ is given by the same way as in {\em Case 3}, except $t$ is replaced by the largest one such that $b_t\in \{i, i+1\}$.
\end{itemize}

For a skew  Young diagram  $\la/\mu$, a tableau $T$ obtained by
filling $\la/\mu$ with letters in $\I$ is called
semistandard  if (1) the letters in each row (resp. column) are
weakly increasing from left to right (resp. from top to bottom), (2)
the letters in $\I_0$ (resp. $\I_1$) are strictly increasing in each
column (resp. row).  We say that the shape of
$T$ is $\la/\mu$. 
We denote by ${SST}(\la/\mu)$ the set of all
semistandard tableaux of shape $\la/\mu$.
It is not difficult to check that $SST(\la)$ is non-empty if and only if $\la\in \cP_{M|N}$.

\begin{ex}{\rm
Suppose that $\e=\e_{3|4}$ with $\I_0=\{1,2,3\}$ and $\I_1=\{4,5,6,7\}$.
$$
\resizebox{.15\hsize}{!}
{\def\lr#1{\multicolumn{1}{|@{\hspace{.75ex}}c@{\hspace{.75ex}}|}{\raisebox{-.04ex}{$#1$}}}\raisebox{-.6ex}
{$\begin{array}{cccccccc}
\cline{1-5}
\lr{\blue {\it 1}} & \lr{\blue {\it 1}} & \lr{\blue {\it 2}} & \lr{\blue {\it 3}} & \lr{6}    \\
\cline{1-5}
\lr{\blue {\it 2}}  & \lr{5} & \lr{7} &  &  \\
\cline{1-3}
\lr{\blue {\it 3}}  & \lr{5} & & & \\
\cline{1-2}
\lr{4}  & \lr{6} & & & \\
\cline{1-2}
\end{array}$}}\quad\quad \in SST((5,3,2,2))
$$
\noindent where the letters in $\I_0$  are in blue (italic).}
\end{ex}

Let $\la\in\cP$ such that $SST(\la)$ is non-empty.
For $T\in SST(\la)$, let $T(i,j)$ denote  the letter of $T$ located in
the $i$-th row from the top and the $j$-th column from the left.
Let $m=\sum_{i\geq 1}\la_i$. Choose an embedding 
\begin{equation}\label{admissible reading}
\iota : SST(\la) \rightarrow \mc{B}^{\otimes m}
\end{equation}
by reading the letters of $T$ in $SST(\la)$
in such a way that $T(i,j)$ should be read before $T(i+1,j)$ and
$T(i,j-1)$ for each $i,j$. Then the image of
$SST(\la)$ under $\iota$ together with $0$ is stable under
$\te_i,\tf_i$ ($i\in \ov I$) and  the induced  $\ov I$-colored oriented graph structure does not depend on the choice of $\iota$ by \cite[Theorem 4.4]{BKK}. 

Now, we suppose that $\epsilon=\epsilon_{M|N}$ in the rest of this subsection.

\begin{prop}\label{prop:crystal of SST}
For $\la\in\cP_{M|N}$, $SST(\la)$ is a connected $\ov I$-colored oriented graph with a unique highest weight element $H_\la$ with highest weight $\la$.
\end{prop}
\pf The first statement is easy to check. The second statement follows from  \cite[Theorem 4.8]{BKK} since the crystal structure on $\mc{B}^{\otimes m}$ and hence on $SST(\la)$ coincides with those in \cite{BKK}.
\qed

\begin{rem}\label{rem:genuine hw/lw vector}
{\rm
We recall that unlike the crystal graph of an integrable highest weight representation of $U_q(\mf{gl}_n)$, there might exist $T\in SST(\la)$ such that ${\rm wt}(T)\neq \la$ but $\te_i T=0$ for all $i\in \ov I$, which is called a fake highest weight vector in \cite{BKK}. 
Due to the existence of fake highest weight vectors,
the proof of the connectedness of $SST(\la)$ becomes much different and non-trivial compared to the case of $U_q(\mf{gl}_n)$. We call the unique highest weight element $H_\la\in SST(\la)$ a genuine highest weight vector. Similarly, there exists a fake lowest weight vector, and we denote by $L_\la$ the genuine lowest weight vector, which is the unique lowest weight vector in $SST(\la)$. }
\end{rem}

For $\la\in\cP_{M|N}$, let
\begin{equation}\label{eq:(L,B)}
\begin{split}
L(\la)&=\sum_{r\geq 0,\, i_1,\ldots,i_r\in \ov I}A_0 \td{x}_{i_1}\cdots\td{x}_{i_r}v_\la, \\
B(\la)&=\{\,\pm \,{\td{x}_{i_1}\cdots\td{x}_{i_r}v_\la}\!\!\! \pmod{q L(\la)}\,|\,r\geq 0, i_1,\ldots,i_r\in \ov I\,\}\setminus\{0\},
\end{split}
\end{equation}
where $v_\la$ is a highest weight vector in $V(\la)$ and $x=e, f$ for each $i_k$. Then we have the following (cf.~\cite[Theorem 5.1]{BKK}).

\begin{thm}\label{thm:crystal base of poly repn}
For $\la\in \cP_{M|N}$, $(L(\la),B(\la))$ is a crystal base of $V(\la)$, and the crystal $B(\la)$ is isomorphic to $SST(\la)$ as an $\ov I$-colored oriented graph.
\end{thm}
\pf By Proposition \ref{prop:irr poly rep}, ${\rm wt}(V(\la))$ is equal to the set of weights of irreducible highest weight $U_q(\mf{gl}_{M|N})$-module with highest weight $\la$. 
Also the $\ov I$-colored oriented graph structure on $SST(\la)$ is equal to the one in \cite{BKK}.
Therefore we may apply the same arguments in \cite[Section 5]{BKK} (together by using Proposition \ref{prop:polarization}) to prove the statements in the theorem. 
\qed

\begin{thm}\label{thm:uniqueness of crystal base}
Let $V$ be a finite-dimensional $\ov\U(\e_{M|N})$-module in $\ov\cO_{\geq 0}$. Suppose that $V$ is completely reducible and it has a crystal base $(L,B)$. Then there exists an isomorphism of $\ov\U(\e_{M|N})$-modules
\begin{equation*}
\psi : V \longrightarrow \bigoplus_{\la\in\cP_{M|N}}V(\la)^{\oplus m_\la},
\end{equation*} 
for some $m_\la\in \Z_{\geq 0}$ such that
\begin{itemize}
\item[(1)] $\psi|_L : L \longrightarrow \bigoplus_{\la\in\cP_{M|N}}L(\la)^{\oplus m_\la}$ is an isomorphism of $A_0$-modules,

\item[(2)] $\ov{\psi}|_{B} : B \longrightarrow \bigsqcup_{\la\in\cP_{M|N}}B(\la)^{\oplus m_\la}$ is an isomorphism of $\ov I$-colored oriented graphs, 
where $\ov{\psi}$ is the map on $L/qL$ induced from $\psi$, and
$B(\la)^{\oplus m_\la}=B(\la)\sqcup\cdots\sqcup B(\la)$ ($m_\la$ times).
\end{itemize}
\end{thm}
\pf Let $\la$ be a maximal weight in ${\rm wt}(V)$ with $\dim V_\la = \ell$. We have a $\Q(q)$-basis of $V_\la$, say $\{\,v_1,\ldots,v_\ell\,\}$ such that $v_i\in L$ and $v_i\in B \pmod{qL}$ for $1\leq i\leq \ell$. 
By Proposition \ref{prop:irr poly rep}, we have $\la\in \cP_{M|N}$ and $\ov\U(\e_{M|N})v_i\cong V(\la)$ for $1\leq i\leq \ell$.
Let $(L(\la)^{(i)},B(\la)^{(i)})$ be the crystal base of $\ov\U(\e_{M|N})v_i$ given in \eqref{eq:(L,B)} generated by $v_i$ for $1\leq i\leq \ell$.

Put $V^{0} = \bigoplus_{1\leq i\leq \ell}\ov\U(\e_{M|N})v_i\cong V(\la)^{\oplus \ell}$. Then  
$L^0:=\bigoplus_{1\leq i\leq \ell}L(\la)^{(i)}$ is a crystal lattice of $V^0$.
On the other hand, $L\cap V^{0}$ is a crystal lattice of $V^0$. By the same arguments as in \cite[Lemma 2.7(iii)]{BKK}, we have $L^0=L\cap V^0$. Since one may regard $L^0/qL^0 \subset L/qL$, we have $B^0:=\bigsqcup_{1\leq i\leq \ell} B(\la)^{(i)}\subset B$. 

Let $V^1$ be a submodule of $V$ such that $V=V^0\oplus V^1$. Put $L^1=L\cap V^1$. Then we have
$L^k/qL^k\subset L/qL$ ($k=0,1$), which is invariant under $\te_i$ and $\tf_i$ for $i\in \ov{I}$. Let $u=\sum_{b\in {\bf B}}c_b b\in L^1/qL^1$ be given, where ${\bf B}$ is a $\Q$-basis of $L/qL$ with $B={\bf B}\cup (-{\bf B})$. If $c_b\neq 0$ for some $b\in B^0$, then $\td{x}_{i_1}\cdots\td{x}_{i_r}u\in (L/qL)_\la\setminus \{0\}=(L^0/qL^0)_\la\setminus\{0\}$ for some $i_1,\ldots, i_r\in \ov{I}$ with $x=e, f$ for each $i_k$, which is a contradiction. Hence $L^1/qL^1$ is spanned by $B^1:=B\setminus B^0$. This implies $L^0/qL^0\oplus L^1/qL^1 = L/qL$ and hence $L=L^0\oplus L^1$. It follows that $(L^1, B^1)$ is a crystal base of $V^1$. The proof completes by using induction on ${\rm dim}V$.
\qed
\vskip 2mm


\begin{cor}\label{cor:LR rule for poly}
We have the following.
\begin{itemize}
\item[(1)] For $\la\in \cP_{M|N}$, $V(\la)$ is a direct summand of $\mc V^{\otimes m}$, where $m=\sum_{i\geq 1}\la_i$, and $V(\la)$ has a polarizable crystal base.

\item[(2)] For $\mu, \nu\in \cP_{M|N}$, $V(\mu)\otimes V(\nu)$ is completely reducible.
The multiplicity of $V(\la)$ in $V(\mu)\otimes V(\nu)$ is equal to the Littlewood-Richardson coefficient associated to $\la$, $\mu$, and $\nu$.
\end{itemize}
\end{cor}
\pf (1) By \cite[Theorem 4.13]{BKK}, we have
\begin{equation*}
\mc{B}^{\otimes m} \cong \bigsqcup_{\substack{\la\in\cP_{M|N}\\ \sum_{i\geq 1}\la_i=m}}
SST(\la)^{\oplus m_\la}, 
\end{equation*}
as an $\ov I$-colored oriented graph, where $m_\la$ is the number of standard tableaux of shape $\la$ with letters in $\{\,1,\ldots, m\,\}$. 
By Corollary \ref{cor:complete reducibility of tensor powers} and Theorem \ref{thm:uniqueness of crystal base}, $V(\la)$ is a submodule of ${\mc V}^{\otimes m}$  
for $\la\in \cP_{M|N}$ with $m=\sum_{i\geq 1}\la_i$. Moreover, by Proposition \ref{prop:polarization} and Theorem \ref{thm:uniqueness of crystal base}, $V(\la)$ has a polarizable crystal base with a polarization induced from the one on $\mc V^{\otimes m}$.

(2) It follows from (1) and Proposition \ref{prop:polarization} that $V(\mu)\otimes V(\nu)$ is completely reducible for $\mu, \nu\in \cP_{M|N}$.
For $\la\in \cP_{M|N}$, the multiplicity of $V(\la)$ in $V(\mu)\otimes V(\nu)$ is equal to the number of connected components in $SST(\mu)\otimes SST(\nu)$ isomorphic to $SST(\la)$, which is shown to be equal to the usual Littlewood-Richardson coefficient asscociated to $\la, \mu, \nu$ \cite[Theorem 4.18]{KK}. 
\qed

\section{$R$ matrix and Kirillov-Reshetikhin modules}\label{sec:KR module}

From now on, we assume that $\e=\e_{M|N}$ for some $M, N\geq 0$ with $M+N=n$. Put $\U=\U(\e_{M|N})$ and $\ov\U=\ov\U(\e_{M|N})$.

\subsection{$R$ matrix}
Let $l, m\in \Z_{\geq 0}$ be given. We suppose that $(l), (m)\in \cP_{M|N}$, which means that $l,m\in \Z_{\ge 0}$ when $M>0$, and $l, m\leq n$ when $M=0$.
For generic $x,y\in \Q(q)$, $\W_l(x)\otimes \W_m(y)$ is an irreducible $\U$-module. 
The case when $M\in\{0, n\}$ is well-known (cf.~\cite{KMN}), and the case for $1\leq M\leq n-1$ is proved in \cite{KOS}.

Consider a non-zero linear map $R\in {\rm End}_{\Q(q)}(\W_l(x)\otimes \W_m(y))$ such that
\begin{equation}\label{eq:R matrix}
\begin{split}
\Delta^{\rm op}(g) \circ R = R \circ  \Delta (g)
\end{split}
\end{equation}
for $g\in \U$, where $\Delta^{\rm op}$ is the opposite coproduct of $\Delta$ in \eqref{eq:comult-1}, that is, $\Delta^{\rm op}(g) = P \circ \Delta(g) \circ P $ and $P(a\otimes b)=b\otimes a$.
Such a map is unique up to scalar multiplication if exists since $\W_l(x)\otimes \W_m(y)$ is irreducible. We call it the {\em $R$ matrix}, and denote it by $R(z)$ since $R$ depends only on $z=x/y$. 

The existence of the $R$ matrix is proved in \cite[Theorem 4.1]{KOS} by reducing a solution of the tetrahedron equation to the Yang-Baxter equation.

\begin{thm}[\cite{KOS}]\label{thm:existence of R-matrix}
For $l, m\in \Z_{\geq 0}$ with $(l), (m)\in \cP_{M|N}$, there exists a non-zero linear map $R\in {\rm End}_{\Q(q)}(\W_l(x)\otimes \W_m(y))$ satisfying \eqref{eq:R matrix} and the Yang-Baxter equation
\begin{equation*}
R_{12}(u)R_{13}(uv)R_{23}(v)=R_{23}(v)R_{13}(uv)R_{12}(u),
\end{equation*}
on $\W_{s_1}(x_1)\otimes \W_{s_2}(x_2)\otimes \W_{s_3}(x_3)$ with $u=x_1/x_2$ and $v=x_2/x_3$ for $(s_1),(s_2),(s_3)\in \cP_{M|N}$ and generic $x_1,x_2,x_3$. Here $R_{ij}(z)$ denotes the map which acts as $R(z)$ on the $i$-th and the $j$-th component and the identity elsewhere.
\end{thm}

Note that $\W_l(x)\cong V((l))$ as a $\ov\U$-module in $\ov\cO_{\geq 0}$, and 
\begin{equation}\label{eq:single row crystal}
\mc{B}_l\cong SST((l))
\end{equation}
as an $\ov I$-colored oriented graph (cf.~\eqref{eq:crystal base of W(r)}), where $|{\bf m}\rangle$ for ${\bf m}=(m_1,\ldots,m_n)$ corresponds to a unique tableau $T$ in $SST((l))$ with $m_i$ the number of occurrences of $i$ in $T$ for $i\in \I$. 

Put 
\begin{equation}\label{eq:H(l,m)}
\begin{gathered}
H(l,m)=\{\,t\,|\,0\leq t\leq \min(l,m), (l+m-t,t)\in \cP_{M|N}\,\}.
\end{gathered}
\end{equation} 
By Corollary \ref{cor:LR rule for poly}, $\W_l(x)\otimes \W_m(y)$ is completely reducible as a $\ov\U$-module, and it has the following multiplicity-free decomposition (see \cite[(6.15)]{KK})
\begin{equation}\label{eq:decomp of W_l and W_m}
\W_l(x)\otimes \W_m(y) \cong 
\bigoplus_{t\in H(l,m)} V((l+m-t,t)).
\end{equation}

For each $t\in H(l,m)$, let $\xi^{l,m}_t$ be the $\ov\U$-highest weight vector of $V((l+m-t,t))$ given in \cite[(6.5), (6.10), (6.14)]{KOS}. More precisely, we have to apply $P$ to $\xi^{l,m}_t$ in \cite{KOS} since our comultiplication is opposite to the one in \cite{KOS}.

We define a $\ov\U$-linear map $\mc P^{l,m}_{t} : \W_l(x)\otimes \W_m(y) \longrightarrow \W_m(y)\otimes \W_l(x)$ given by $\mc P^{l,m}_{t}(\xi^{l,m}_{t'})=\delta_{tt'}\xi^{m,l}_t$ for $t'\in H(l,m)$. Then the $\U$-linear map $PR(z)$ has the spectral decomposition 
\begin{equation}\label{eq:spectral decomp}
PR(z)= \sum_{t\in H(l,m)}\rho_t(z) \mc P^{l,m}_t,
\end{equation}
for some $\rho_t(z)\in \Q(q)$. An explicit form of $PR(z)$ is given as follows \cite[(6.10), (6.13), (6.16)]{KOS};

\begin{equation}\label{eq:spectral M=0}
PR(z) =  
\sum_{t=\max\{l+m-n,0\}}^{\min(l,m)}\left(\prod_{i=t+1}^{\min(l,m)}\frac{1-q^{l+m-2i+2}z}{z-q^{l+m-2i+2}}\right)\mc P^{l,m}_t \ \quad (M=0),
\end{equation}

\begin{equation}\label{eq:spectral M=1}
PR(z)=
\sum_{t=0}^{\min(l,m,n-1)}\left(\prod_{i=1}^{t}\frac{z-q^{l+m-2i+2}}{1-q^{l+m-2i+2}z}\right)\mc P^{l,m}_t\hskip 2cm (M=1),
\end{equation}

\begin{equation}\label{eq:spectral M>1}
PR(z)=
\sum_{t=0}^{\min(l,m)}
\left(\prod_{i=1}^{t}\frac{z-q^{l+m-2i+2}}{1-q^{l+m-2i+2}z}\right)\mc P^{l,m}_t
\hskip 2cm (2\leq M\leq n),
\end{equation}
where we assume that $\rho_{\min(l,m)}(z)=1$ in \eqref{eq:spectral M=0} and $\rho_0(z)=1$ in \eqref{eq:spectral M=1} and \eqref{eq:spectral M>1}.
We remark that $\rho_t(z)$ above is equal to $\rho_t(z^{-1})$ in \cite{KOS} by our convention of comultiplication.

\subsection{Fusion construction of KR modules}\label{subsec:fusion construction}
Let us apply fusion construction following \cite{KMN2}. 
We assume that $1\leq M \leq n-1$ since the results when $M\in \{0, n\}$ are well-known (cf.~ \cite{KMN2}). 

Fix $s\geq 1$ and put $V_x=\mc{W}_s(x)$ for $x\in \Q(q)$.
For  $R(z)\in {\rm End}_{\Q(q)}(V_x\otimes V_y)$ in \eqref{eq:spectral M=1} and \eqref{eq:spectral M>1}, we take a normailzation
\begin{equation*}
\check{R}(z) = \left(\prod_{i=1}^s \frac{1-q^{2s-2i+2}z}{z-q^{2s-2i+2}}\right) PR(z).
\end{equation*}
Since $(s^2)\not\in \cP_{M|N}$ if and only if $M=1$ and $s>n-1$, we have
\begin{equation}\label{eq:normalized R}
\check{R}(z)=
\begin{cases}
\sum_{t=0}^{n-1}\left(\prod_{i=t+1}^{s}\dfrac{1-q^{2s-2i+2}z}{z-q^{2s-2i+2}}\right)\mc P^{s,s}_t,& \text{if $(s^2)\not\in \cP_{M|N}$},\\
\mc P^{s,s}_s + \sum_{t=0}^{s-1}\left(\prod_{i=t+1}^{s}\dfrac{1-q^{2s-2i+2}z}{z-q^{2s-2i+2}}\right)\mc P^{s,s}_t,& \text{if $(s^2)\in \cP_{M|N}$}. 
\end{cases}
\end{equation}
For $r\geq 2$, let $\mf{S}_r$ denote the group of permutations on $r$ letters generated by $s_i=(i\ i+1)$ for $1\leq i\leq r-1$.
By Theorem \ref{thm:existence of R-matrix}, we have $\U$-linear maps 
\begin{equation*}
\check{R}_w(x_1,\ldots,x_r) : 
V_{x_1}\otimes \cdots\otimes V_{x_r} \longrightarrow V_{x_{w(1)}}\otimes \cdots\otimes V_{x_{w(r)}}
\end{equation*}
for $w\in \mf{S}_r$ and generic $x_1,\ldots, x_r$ satisfying the following:
\begin{align*}
&\check{R}_1(x_1,\ldots,x_r) = {\rm id}_{V_{x_1}\otimes \cdots\otimes V_{x_r}},\\
&\check{R}_{s_i}(x_1,\ldots,x_r) = \left(\otimes_{j<i}{\rm id}_{V_{x_j}}\right)\otimes \check{R}(x_i/x_{i+1}) \otimes \left(\otimes_{j>i+1}{\rm id}_{V_{x_j}}\right),\\
&\check{R}_{ww'}(x_1,\ldots,x_r) = \check{R}_{w'}(x_{w(1)},\ldots,x_{w(r)})\check{R}_{w}(x_1,\ldots,x_r),
\end{align*}
for $w, w'\in \mf{S}_r$ with $\ell(ww')=\ell(w)+\ell(w')$, where $\ell(w)$ denotes the length of $w$.

Let $w_0$ denote the longest element in $\mf{S}_r$. By \eqref{eq:normalized R}, $\check{R}_{w_0}(x_1,\ldots,x_r)$ does not have a pole at $q^{-2k}$ for $k\in \Z_{\geq 0}$ as a function in $x_1,\ldots,x_r$. Hence we have a $\U$-linear map
\begin{equation*}
\begin{split}
\check{R}_r := \check{R}_{w_0}(q^{1-r},q^{3-r},\cdots, q^{r-1}) : V_{q^{1-r}}\otimes \cdots\otimes V_{q^{r-1}} \longrightarrow 
V_{q^{r-1}}\otimes \cdots\otimes V_{q^{1-r}}.
\end{split}
\end{equation*}
We define a $\U$-module
\begin{equation*}
\W_s^{(r)} = {\rm Im}\check{R}_r,
\end{equation*}
which we call a {\em Kirillov-Reshetikhin module} corresponding to $(r,s)$.

\begin{prop}\label{prop:W_s^{(r)}}
For $r, s\geq 1$, the $\U$-module $\W_s^{(r)}$ is non-zero if and only if $(s^r)=(\underbrace{s,\cdots,s}_{r\ {\rm times}})\in \cP_{M|N}$. In this case, we have $\W_s^{(r)}\cong V((s^r))$ as a $\ov\U$-module.
\end{prop}
\pf The proof is similar to the case when $M=0$ in \cite[Section 3]{KMN2}. But we give a self-contained proof for the reader's convenience.
We have $V_x \cong V_1\cong V((s))$ as a $\ov\U$-module for $x\in \Q(q)$ with a crystal base $(L((s)),B((s)))$. 
By Corollary \ref{cor:LR rule for poly}(2) and \cite[(6.15)]{KK}, we have as a $\ov\U$-module
\begin{equation}\label{eq:decomp of tensor power of V(r)}
\begin{split}
V_{q^{1-r}}\otimes \cdots\otimes V_{q^{r-1}} &\cong 
V_{q^{r-1}}\otimes \cdots\otimes V_{q^{1-r}} \cong V((s))^{\otimes r}  \cong   \bigoplus_{\la}V(\la)^{\oplus m_\la},
\end{split}
\end{equation}
where the sum is over $\la\in \cP_{M|N}$ such that $\sum_{i\geq 1}\la_i=rs$ and $\la_i\geq s$ for some $i$ and $m_\la\in \Z_{\geq 0}$. Here we have $m_{(s^r)}\leq 1$, and $m_{(s^r)}=1$ if and only if $(s^r)\in \cP_{M|N}$.

For each $1\leq i\leq r-1$, we have
\begin{equation*}
\check{R}_r = R_{s_i}(\cdots, \underbrace{q^{r-2i-1}}_{i},\underbrace{q^{r-2i+1}}_{i+1},\cdots) R_{w_0s_i}(q^{1-r},\cdots, q^{r-1}),
\end{equation*}
which implies that as a $\ov\U$-module
\begin{equation}\label{eq:image of R_s}
\W_s^{(r)} \subset V((s))^{\otimes i-1}\otimes {\rm Im}\check{R}_2 \otimes V((s))^{\otimes r-i-1}\subset V((s))^{\otimes r}.
\end{equation}
By \eqref{eq:normalized R}, we have as a $\ov\U$-module ${\rm Im}\check{R}_2 = V((s^2))$ if $(s^2)\in \cP_{M|N}$, and ${\rm Im}\check{R}_2 = 0$ otherwise.

Suppose that $\W_s^{(r)}$ is non-zero.
By Theorem \ref{thm:uniqueness of crystal base}, $\W^{(r)}_s$ as a $\ov\U$-module has a crystal base $(L',B')$, where $L'$ is a direct summand of $L((s))^{\otimes r}$ and $B'$ is a subset of $B((s))^{\otimes r}$. 
Let $b_1\otimes \cdots \otimes b_r\in B'$ be given. 
Since ${\rm Im}\check{R}_2 = V((s^2))$ by \eqref{eq:image of R_s}, 
$C(b_i\otimes b_{i+1})$, the connected component of $b_i\otimes b_{i+1}$ in $B((s))^{\otimes 2}$ with respect to $\te_i$ and $\tf_i$ for $i\in \ov I$, is isomorphic to $B((s^2))\cong SST((s^2))$ as an $\ov I$-colored oriented graph. 
This is equivalent to saying that $b_i\otimes b_{i+1}$ corresponds to a tableau in $SST((s^2))$, whose first and second row correspond to $b_i$ and $b_{i+1}$, respectively.
Hence $C(b_1\otimes \cdots \otimes b_r)\subset B((s))^{\otimes r}$ is isomorphic to $B((s^r))\cong SST((s^r))$. It follows that $(s^r)\in \cP_{M|N}$.

Conversely, suppose that $(s^r)\in \cP_{M|N}$. 
By \eqref{eq:decomp of tensor power of V(r)}, 
there exists a unique highest weight vector $u_+$ with highest weight $\La_{(s^r)}$ (see \eqref{eq:highest weight correspondence}) in $V_1^{\otimes r}$ up to scalar multiplication. When $r=2$, we have $\check{R}(q^{-2k})(u_+)\ne0$ for any integer
$k\ge1$ from \eqref{eq:normalized R}. In view of the definition of $\check{R}_r$, we 
also have $\check{R}_r(u_+)\ne0$ for any $r$, and hence $\W_s^{(r)}$ is non-zero. 

Since the crystal base of $\W_s^{(r)}$ as a $\ov{\U}$-module is contained in $B((s^r))$ and $V((s^r))$ is irreducible, we have $\W_s^{(r)}=V((s^r))$.
\qed\vskip 2mm


Let $r,s\geq 1$ such that $(s^r) \in \cP_{M|N}$.
Let $u_+\in \W_s^{(r)}$ be a $\ov{\U}$-highest weight vector corresponding to $H_{(s^r)}\in SST((s^r))$ and 
let $u_-\in \W_s^{(r)}$ be a $\ov{\U}$-lowest weight vector corresponding to $L_{(s^r)}\in SST((s^r))$ (see Remark \ref{rem:genuine hw/lw vector}). Put
\begin{equation}\label{eq:extremal vector}
u_0 =
\begin{cases}
u_+ & \text{if $r\leq M$},\\
u_- & \text{if $r > M$},\\
\end{cases}
\end{equation}
where we assume $u_0\in L((s))^{\otimes r}$ and $u_0 \in B((s))^{\otimes r} \pmod {q L((s))^{\otimes r}}$.  

\begin{ex}{\rm
Suppose that $n=7$ and $\e=\e_{3|4}$, where $\I_0=\{1,2,3\}$ and $\I_1=\{4,5,6,7\}$.
Then
$$u_0=u_+\equiv
\resizebox{.15\hsize}{!}
{\def\lr#1{\multicolumn{1}{|@{\hspace{.75ex}}c@{\hspace{.75ex}}|}{\raisebox{-.04ex}{$#1$}}}\raisebox{-.6ex}
{$\begin{array}{ccccc}
\cline{1-5}
\lr{1} & \lr{1} &  \lr{1} & \lr{1} & \lr{1} \\
\cline{1-5}
\lr{2} & \lr{2} &  \lr{2} & \lr{2} & \lr{2} \\
\cline{1-5}
\lr{3} & \lr{3} &  \lr{3} & \lr{3} & \lr{3} \\
\cline{1-5}
\end{array}$}}\in SST((5^3))\quad \quad \quad  
u_0=u_-\equiv
\resizebox{.09\hsize}{!}
{\def\lr#1{\multicolumn{1}{|@{\hspace{.75ex}}c@{\hspace{.75ex}}|}{\raisebox{-.04ex}{$#1$}}}\raisebox{-.6ex}
{$\begin{array}{ccc}
\cline{1-3}
\lr{5} & \lr{6} & \lr{7}  \\
\cline{1-3}
\lr{5} & \lr{6} & \lr{7}  \\
\cline{1-3}
\lr{5} & \lr{6} & \lr{7}  \\
\cline{1-3}
\lr{5} & \lr{6} & \lr{7}  \\
\cline{1-3}
\end{array}$}}\in SST((3^4))\quad .
$$\vskip 2mm

}
\end{ex}

\begin{prop}\label{prop:polarization on KR}
There exists a polarization $(\ , \ )_r$ on $\W_s^{(r)}$ such that $(u_0,u_0)_r=1$.
\end{prop}
\pf By Proposition \ref{prop:polarization on W(r)}, the polarization on $\W_s=V_1$ in \eqref{eq:polarization on W(r)} yields a pairing on $V_x\otimes V_{x^{-1}}$ for $x\in \Q(q)$ satisfying the condition \eqref{eq:polarization}.
Hence we have a pairing on 
$(V_{x_1}\otimes \cdots \otimes V_{x_r})\otimes (V_{x_1^{-1}}\otimes \cdots \otimes V_{x_r^{-1}})$ for $x_1,\ldots,x_r\in \Q(q)$ given by
\begin{equation}\label{eq:pairing-2}
(u,v) = (u_1,v_1)\cdots (u_r,v_r),
\end{equation}
for $u=u_1\otimes \cdots \otimes u_r\in V_{x_1}\otimes \cdots \otimes V_{x_r}$ and $v=v_1\otimes \cdots \otimes v_r\in V_{x_1^{-1}}\otimes \cdots \otimes V_{x_r^{-1}}$, which  also satisfies \eqref{eq:polarization}. By putting $(x_1,x_2,\ldots,x_r) = (q^{1-r},q^{3-r},\cdots,q^{r-1})$ in \eqref{eq:pairing-2}, we have a pairing on 
$(V_{q^{1-r}}\otimes \cdots \otimes V_{q^{r-1}})\otimes (V_{q^{r-1}}\otimes \cdots \otimes V_{q^{1-r}})$.

Let $\widetilde{V}=V_{q^{r-1}}\otimes \cdots \otimes V_{q^{1-r}}$.
By Proposition \ref{prop:W_s^{(r)}}, $\W_s^{(r)}$ is a $\U$-submodule of $\widetilde{V}$ generated by $u_0$.
Note that $\widetilde{V}\cong V((s))^{\otimes r}$ as a $\ov\U$-module.
The crystal base $(L((s)),B((s)))$ of $\W_s=V((s))$ as a $\ov\U$-module is polarizable with respect to \eqref{eq:polarization on W(r)} by Proposition \ref{prop:polarization on W(r)}, and hence $(L((s))^{\otimes r},B((s))^{\otimes r})$ is a polarizable crystal base of $\widetilde{V}$ with respect to \eqref{eq:pairing-2} by Proposition \ref{prop:polarization}. 
So, we have $(u_0,\check{R}_r(u_0)) = c_r$ for some $c_r \in 1+qA_0$. Now, for $u,v\in \W_s^{(r)}$, we define 
$(u,v)_r = c_r^{-1} (u,\check{R}_r(v))$, where $(u,\check{R}_r(v))$ is given in \eqref{eq:pairing-2}. One can check that $(\ ,\ )_r$ is symmetric since $\W_s^{(r)}$ is irreducible $\ov\U$-module generated by $u_0$ (cf.~\cite[Lemma 3.4.2]{KMN2}). Hence $(\ ,\ )_r$ is a polarization on $\W_s^{(r)}$.
\qed\vskip 2mm

\subsection{Crystal base of KR modules}
We assume as in Section \ref{subsec:fusion construction} that $1\leq M\leq n-1$.
Let $r, s\geq 1$ such that $(s^r)\in \cP_{M|N}$. 

Let $A_{\Z} = \{\,f(q)/g(q)\in\Q(q)\,|\, f(q), g(q)\in \Z[q],\ g(0)=1\,\}$ and  
$K_{\Z} = A_\Z[q^{-1}]$. 
Note that $A_0\cap K_\Z=A_\Z$. 
We define $\U_{K_\Z}$ to be the $K_\Z$-subalgebra of $\U$ generated by 
$e_i^{(a)}$, $f_i^{(a)}$, and $q^h$ for $i\in I$, $a\geq 0$ and $h\in P^\vee$.

Put $W = \W_s^{(r)}$ and let $(\ ,\ )_r$ be the polarization on $W$ in Proposition \ref{prop:polarization on KR}. Let $W_{K_\Z} = \U_{K_\Z} u_0$, where $u_0$ is as in \eqref{eq:extremal vector}. 

\begin{lem}
We have $(W_{K_\Z},W_{K_\Z})_r\subset K_\Z$.
\end{lem}
\pf Since 
$(W_{K_\Z},W_{K_\Z})_r = (\U_{K_\Z} u_0, \U_{K_\Z} u_0)_r =( u_0, \U_{K_\Z} u_0)_r$,
and $(u_0,u_0)_r=1$, it is enough to show that $\U_{K_\Z} u_0 \cap \Q(q) u_0 = K_\Z u_0$.

{\em Case 1}. Suppose that $r\leq M$, and $u_0=u_+$. 
Note that $u_0$ is a highest weight vector with respect to the subalgebra of $\ov\U$ generated by $e_i$, $f_i$, and $q^{\delta_j^\vee}$ for $i=1,\ldots, r-1$ and $j=1,\ldots,r$, which is isomorphic to $U_q(\mf{gl}_r)$. We may take $u_+$ as the highest weight vector in the upper global crystal basis of the polynomial representation of $U_q(\mf{gl}_r)$ with highest weight $s(\delta_1+\cdots+\delta_r)$ constructed in \cite{Br}. (We take a basis element $L_{R(A)}$ in \cite[Theorem 26]{Br} corresponding to a highest weight vector. Here we replace $q^{-1}$ appearing in $L_{R(A)}$ with $q$ which is due to the difference of comultiplication.) So $u_+$ lies in the tensor product of $q$-deformed symmetric powers of the natural representation of $U_q(\mf{gl}_r)$. 
Moreover, we have
\begin{equation}\label{eq:integral form of u_s}
u_0 =u_+= \sum_{\substack{{\bf m}_1,\cdots,{\bf m}_r \in \Z_{\geq 0}^n \\ |{\bf m}_i|=s}} c_{{\bf m}_1,\cdots,{\bf m}_r}(q)|{\bf m}_1\rangle \otimes \cdots \otimes |{\bf m}_r\rangle,
\end{equation}
for some $c_{{\bf m}_1,\cdots,{\bf m}_r}(q)\in \Z[q]$ such that
\begin{equation}\label{eq:integral form of u_s-2}
c_{{\bf m}_1,\cdots,{\bf m}_r}(0)=
\begin{cases}
1, & \text{if ${\bf m}_i=s\be_i$ for $1\leq i\leq r$},\\
0, & \text{otherwise}.
\end{cases}
\end{equation}
Let $V_{K_\Z}=\bigoplus_{|{\bf m}|=s}K_\Z |{\bf m}\rangle$. Since $\U_{K_\Z} V_{K_\Z}\subset V_{K_\Z}$, we have $\U_{K_\Z} u_0 \subset (V_{K_\Z})^{\otimes r}$ by \eqref{eq:integral form of u_s}.
Suppose that $c(q)u_0\in \U_{K_\Z} u_0 \cap \Q(q) u_0$. We may assume that $c(q)=c_1(q)/c_2(q)$ where $c_1(q), c_2(q)\in \Z[q]$ and they are relatively prime. 
Let $c_0(q)=c_{s\be_1,\ldots,s\be_r}(q)$.
Then we have $c(q)\in K_{\Z}$ since $c_0(q)c(q)\in K_\Z$ and $c_0(0)=1$.

{\em Case 2}. Suppose that $r > M$, and $u_0=u_-$.
Clearly, we have 
\begin{equation*}
u_0=u_-=|\be_{s-n+1}+\dots+\be_n\rangle^{\otimes r}. 
\end{equation*}
By the same argument as in {\em Case 1}, we have $\U_{K_\Z} u_0 \cap \Q(q) u_0 = K_\Z u_0$.
\qed


\vskip 2mm

Let $\Q(q)_+ = \bigsqcup_{l\in\Z,\, c>0} q^l(c + q A_0) \subset \Q(q)$. For $f(q), g(q)\in \Q(q)$, we define $f(q)\geq g(q)$ if $f(q)-g(q) \in \Q(q)_+$. Let $(\ ,\ )$ be a $\Q(q)$-valued symmetric bilinear form on a $\Q(q)$-space $V$. We say that $(\ , \ )$ is positive definite if $(v,v)>0$ for any non-zero $v\in V$.

\begin{lem}\label{lem:positive definiteness}
The polarization $(\ , \ )_r$ is positive definite. 
\end{lem}
\pf It follows from the proof of Proposition \ref{prop:polarization on KR} 
and Theorem \ref{thm:uniqueness of crystal base} that $(L((s^r)),B((s^r)))$ is a crystal base of $W=V((s^r))$ as a $\ov\U$-module, and it is polarizable with respect to $(\ ,\ )_r$. Hence $(\ , \ )_r$ is positive definite by \cite[Lemma 2.2.2]{KMN2}.  
\qed

\begin{lem}\label{lem:norm}
Let $u, u'\in W$ such that ${\rm wt}(u)={\rm wt}(u')=\la$ and $e_iu=e_iu'=0$ for some $\la\in P$ and $i\in I_{\rm even}$.
For $a\geq 0$, we have
\begin{equation*}
(f_i^{(a)}u, f_i^{(a)}u')_r = 
q^{\si a(a-\langle\la,\alpha^\vee_i\rangle)} 
\begin{bmatrix} \langle\la,\alpha^\vee_i\rangle \\ a \end{bmatrix} (u, u')_r.
\end{equation*}
\end{lem}
\pf
We have for $a\geq 0$
\begin{equation}\label{eq:ef commutation}
e_i^{(a)}f_i^{(a)} = \sum_{l=0}^a f_i^{(a-l)}e_i^{(a-l)}\stirlingii{k_i}{l},
\end{equation}
where
\begin{equation*}
\stirlingii{k_i}{l} = \frac{\{k_i\}\{q^{-\si}k_i\}\cdots \{q^{(1-l)\si}k_i\}}{[l]!},
\end{equation*} 
with $\{x\}=\frac{x-x^{-1}}{q-q^{-1}}$ for $x\in q^sk_i$ $(s\in \Z)$.
Note that \eqref{eq:ef commutation} is well-known when $(\e_i,\e_{i+1})=(0,0)$, and the case when $(\e_i,\e_{i+1})=(1,1)$ is obtained by applying $\xi$ in \eqref{eq:sl_2 -} to it.

For $a\geq 0$, we have
\begin{equation*}
\begin{split}
(f_i^{(a)}u, f_i^{(a)}u')_r& = (\eta(f_i^{(a)})f_i^{(a)}u, u')_r= 
((-1)^{a\e_i} q^{\si a^2} k_i^{-a} e_i^{(a)}f_i^{(a)}u, u')_r \\
&= ((-1)^{a\e_i} q^{\si a^2} k_i^{-a} \stirlingii{k_i}{a} u, u')_r \\
& = (-1)^{a\e_i} (-1)^{a(\la_i\e_i +\la_{i+1}\e_{i+1})} q^{\si a(a-\langle\la,\alpha^\vee_i\rangle)}  (\stirlingii{k_i}{a} u, u')_r \\
&= q^{\si a(a-\langle\la,\alpha^\vee_i\rangle)} 
\begin{bmatrix} \langle\la,\alpha^\vee_i\rangle \\ a \end{bmatrix} (u, u')_r, 
\end{split}
\end{equation*}
since 
\begin{equation*}
\stirlingii{k_i}{a} u = 
(-1)^{a\e_i} (-1)^{a(\la_i\e_i +\la_{i+1}\e_{i+1})}
\begin{bmatrix} \langle\la,\alpha^\vee_i\rangle \\ a \end{bmatrix}u.
\end{equation*}

\qed

\begin{lem}\label{lem:norm inequality}
For $u\in W$ and $i\in I$, we have
\begin{equation*}
\begin{split}
&(\te_i u, \te_iu)_r \leq (1+q)(u,u)_r, \quad (\tf_i u, \tf_iu)_r \leq (1+q)(u,u)_r.
\end{split}
\end{equation*}
\end{lem}
\pf Suppose that $u\in W_\la$ for some $\la\in P$. 
We prove only $(\tf_i u, \tf_iu)_r \leq (1+q)(u,u)_r$ 
since the proof of $(\te_i u, \te_iu)_r \leq (1+q)(u,u)_r$ is similar.

{\em Case 1}. $i\in I_{\rm even}$ and $(\e_i,\e_{i+1})=(1,1)$.
Let $u=\sum_{k\geq 0}f_i^{(k)}u_k$, where $e_i u_k=0$ for $k\geq 0$.  
By Lemma \ref{lem:norm}, we have
\begin{equation}\label{eq:norm-1}
(f_i^{(k)}u_k, f_i^{(k)}u_k)_r = 
q^{k(\langle\la,\alpha^\vee_i\rangle+k)} 
\begin{bmatrix} \langle\la,\alpha^\vee_i\rangle +2k \\ k \end{bmatrix} (u_k, u_k)_r.
\end{equation}
Then we may apply \cite[Proposition 2.4.4]{KMN2} to have $(\tf_i u, \tf_iu)_r \leq (1+q)(u,u)_r$.

{\em Case 2}. $i\in I_{\rm even}$ and $(\e_i,\e_{i+1})=(0,0)$. 
Let $u=\sum_{k\geq 0}f_i^{(k)}u_k$, where $e_i u_k=0$ for $k\geq 0$.  
If we put $w_k = q^{-k(l_k-k)}u_k$ where $l_k =\langle \la + k\alpha_i,\alpha_i^\vee \rangle$, then 
\begin{equation*}
\begin{split}
u&=\sum_{k\geq 0}f_i^{(k)}u_k = \sum_{k\geq 0}q^{k(l_k-k)} f_i^{(k)}w_k=\sum_{k\geq 0} \tf_i^{k}w_k,\\
\tf_i u&=\sum_{k\geq 0}q^{(k+1)(l_k-k-1)} f_i^{(k+1)}w_k=\sum_{k\geq 0} \tf_i^{k+1}w_k.
\end{split}
\end{equation*}
By Lemma \ref{lem:norm}, we have
\begin{equation}\label{eq:norm-2}
(\tf_i^{k}w_k, \tf_i^{k}w_k)_r = 
q^{k(\langle\la,\alpha^\vee_i\rangle+k)} 
\begin{bmatrix} \langle\la,\alpha^\vee_i\rangle +2k \\ k \end{bmatrix} (w_k, w_k)_r.
\end{equation}
and apply the same argument as in {\em Case 1}.

{\em Case 3}. $i\in I_{\rm odd}$ and $(\e_i,\e_{i+1})=(1,0)$.
Let $u=u_0 + f_i u_1$, where $e_i u_k=0$ for $k= 0,1$. 
We may assume that  $\langle \la,\alpha_i^\vee \rangle =l>0$.
Then we have $(u,u)_r = (u_0,u_0)_r + (f_iu_1,f_iu_1)_r$ and 
\begin{equation}\label{eq:norm-3}
\begin{split}
(\tf_i u,\tf_i u)_r & = (f_i u_0,f_i u_0)_r =  (\eta(f_i)f_iu_0,u_0)_r = (q_ik_i^{-1}e_if_iu_0,u_0)_r \\
& = (q_ik_i^{-1} \frac{k_i-k_i^{-1}}{q-q^{-1}}u_0,u_0)_r = (\frac{1-k_i^{-2}}{1-q^2}u_0,u_0)_r \\
& = \frac{1-q^{2l}}{1-q^2} (u_0,u_0)_r \leq (1+q)(u_0,u_0)_r\leq (1+q)(u,u)_r. 
\end{split}
\end{equation}

{\em Case 4}. $i\in I_{\rm odd}$ and $(\e_i,\e_{i+1})=(0,1)$. 
Let $u=u_0 + e_i u_1$, where $f_i u_k=0$ for $k= 0,1$. 
We may assume that $\langle \la,\alpha_i^\vee \rangle =l>0$. Put $v_1=e_iu_1$.
Then we have $(u,u)_r = (u_0,u_0)_r + (v_1,v_1)_r$ and 
\begin{equation}\label{eq:norm-4}
\begin{split}
(\tf_i u,\tf_i u)_r & = (\eta(e_i)v_1,\eta(e_i)v_1)_r =  (e_i\eta(e_i)v_1,v_1)_r = (e_iq_i^{-1}f_ik_iv_1,v_1)_r \\
& =  (q_i^{-1}k_ie_if_iv_1,v_1)_r  = (q_i^{-1}k_i \frac{k_i-k_i^{-1}}{q-q^{-1}}v_1,v_1)_r = (\frac{k_i^{2}-1}{q^2-1}v_1,v_1)_r \\
& = \frac{q^{2l}-1}{q^2-1} (v_1,v_1)_r \leq (1+q)(v_1,v_1)_r\leq (1+q)(u,u)_r. 
\end{split}
\end{equation}
\qed

Let 
\begin{equation}\label{eq:crystal lattice of KR}
\begin{split}
L^{r,s}&=\{\,u\in W\,|\, (u,u)_r\in A_0\,\},\\
B^{r,s}&=\{\,b\in W_{K_\Z}\cap L^{r,s} /qW_{K_\Z}\cap L^{r,s}\,|\, (b,b)_0=1\,\},
\end{split}
\end{equation}
where $(\ ,\ )_0$ denotes the bilinear form on $L^{r,s} /q L^{r,s}$ induced from $(\ , \ )_r$.
We assume that $(L((s^r)),B((s^r)))$ is the crystal base of $W=V((s^r))$ as a $\ov\U$-module such that $u_0\in B((s^r)) \pmod{qL((s^r))}$.

\begin{lem}\label{lem:crystal lattice of KR}
$L^{r,s}$ is a crystal lattice of $W$, and $L^{r,s} = L((s^r))$.
\end{lem}
\pf It follows immediately from Lemma \ref{lem:norm inequality} that $L^{r,s}$ is a crystal lattice of $W$. 
Also, we have shown in the proof of Lemma \ref{lem:positive definiteness} that $(\ , \ )_0$ is positive definite on $L((s^r))/qL((s^r))$. Hence $L^{r,s} = L((s^r))$ by \cite[Lemma 2.2.2]{KMN2}.
\qed

\begin{lem}\label{lem:adjont at q=0}
For $u, v\in L^{r,s}$ and $i\in I$,
$(\te_i u, v)_0 = (u, \tf_i v)_0$.
\end{lem}
\pf {\em Case 1}. $i\in I_{\rm even}$. 
We may assume that $u=\tf_i^{k+1} u_0$ and $v=\tf_i^k v_0$ for some $k\geq 0$ and $u_0$, $v_0$ such that $e_i u_0=e_i v_0=0$. 
Let ${\rm wt}(u_0)=\la$.
By \eqref{eq:norm-1} and \eqref{eq:norm-2}, we have
\begin{equation*}
(\tf_i^{k+1}u_0, \tf_i^{k+1}v_0)_r=
\frac{1-q^{2(\langle \la,\alpha_i^\vee \rangle-k)}}{1-q^{2(k+1)}}(\tf_i^{k}u_0, \tf_i^{k}v_0)_r,
\end{equation*}
which implies that $(\tf_i^{k}u_0, \tf_i^{k}v_0)_0 = (\tf_i^{k+1}u_0, \tf_i^{k+1}v_0)_0$, and hence $(\te_i u, v)_0 = (u, \tf_i v)_0$.

{\em Case 2}. $i\in I_{\rm odd}$. Assume 
$u=\tf_i u_0$ and $v= v_0$ for some $u_0$ and $v_0$ such that $e_i u_0=e_i v_0=0$. 
By \eqref{eq:norm-3} and \eqref{eq:norm-4}, we have
$(u_0, v_0)_0=(\tf_iu_0, \tf_iv_0)_0$, and hence $(\te_i u, v)_0 = (u, \tf_i v)_0$.
\qed\vskip 2mm

Now we have the following, which is one of the main results in this paper.

\begin{thm}\label{thm:crystal base of KR}
$(L^{r,s}, B^{r,s})$ is a crystal base of $\W_s^{(r)}$, and $B^{r,s} = B((s^r))$.
\end{thm}
\pf 
Let $b\in B^{r,s}$ such that $\te_i b\neq 0$ for some $i\in I$.  
Let $u\in W_{K_\Z}\cap L^{r,s}$ such that $u\equiv  b \pmod{qL^{r,s}}$. 
Since $W_{K_\Z}\cap L^{r,s}$ is invariant under $\te_i$ and $\tf_i$ for $i\in I$, $\te_i u \in W_{K_\Z}\cap L^{r,s}$. By Lemma \ref{lem:norm inequality}, we have $(\te_i b,\te_i b)_0=1$ and hence $\te_ib\in B^{r,s}$. Similarly, we have  $\tf_i b \in B^{r,s}\cup \{0\}$. Therefore, $B^{r,s}\cup\{0\}$ is invariant under $\te_i$ and $\tf_i$ for $i\in I$.

For $b\in B^{r,s}$ and $i\in I$, if $\te_i b\neq 0$, then $(\tf_i\te_i b, b)_0=(\te_i b,\te_i b)_0=1$ by Lemma \ref{lem:adjont at q=0}. In particular, $\tf_i\te_ib\neq 0$ and $(\tf_i\te_ib-b,\tf_i\te_ib-b)_0=0$, which implies that $\tf_i\te_ib=b$. The proof of $\te_i\tf_ib=b$ when $\tf_ib\neq 0$ is the same.

Let $b_0\in W_{K_\Z}\cap L^{r,s}/qW_{K_\Z}\cap L^{r,s}$ such that $u_0\equiv b_0 \pmod{q L^{r,s}}$. Then $b_0\in B^{r,s}$ since $u_0\in W_{K_\Z}\cap L^{r,s}$ by Proposition \ref{prop:polarization on KR}. 
Let $C_0$ be the connected component of $b_0$ in $B^{r,s}$ with respect to $\te_i$ and $\tf_i$ for $i\in \ov I$.
By Lemma \ref{lem:crystal lattice of KR} and Proposition \ref{prop:crystal of SST}, $C_0=B((s^r))$. So $B^{r,s}$ spans $L^{r,s}/qL^{r,s}$ and accordingly it is a pseudo-basis of $L^{r,s}/qL^{r,s}$ by \cite[Lemma 2.3.2]{KMN2}. In particular $B^{r,s}=B((s^r))$.  Therefore, $(L^{r,s},B^{r,s})$ is a crystal base of $W$.
\qed

\section{Combinatorial description of Kirillov-Reshetikhin crystals}\label{sec:KR crystal}

We assume that $\U=\U(\e_{M|N})$ and $\ov\U=\ov\U(\e_{M|N})$ with $M+N=n$ and $1\leq M\leq n-1$.

\subsection{Subalgebra $\ov\U^\sigma$}
Let 
$$
\e'_{N|M}=(\underbrace{1,\cdots,1}_{N},\underbrace{0,\cdots,0}_{M}).
$$
Let $\sigma : \U(\e'_{N|M})\longrightarrow \U(\e_{M|N})$
be the $\Q(q)$-algebra isomorphism given by $\sigma(e_i)=e_{i+M}$, $\sigma(f_i)=f_{i+M}$, and $\sigma(q^{\de_j^\vee})=q^{\de_{j+M}^\vee}$ for $i\in I$ and $j\in \I$.
Here the subscripts $i,j$ are understood to be modulo $n$. 

Put $\ov I^\sigma:=I\setminus\{M\}$.
Let $\ov\U^\sigma=\sigma\left(\ov{\U}(\e'_{N|M})\right)$ be the subalgebra of $\U$ generated by $q^h$, $e_i, f_i$ for $h\in P^\vee$ and $i\in \ov I^\sigma$. 
The results in Sections \ref{subsec:Yamane's quantum group} and  \ref{subsec:crystal of poly repn} also hold in case of $\ov\U^\sigma$. Let us briefly summarize them without proof since they can be proved in almost the same way as in case of $\ov\U(\e_{M|N})$. 


Let $\ov{\mc{O}}_{\geq 0}^\sigma$ be the category of $\ov\U^\sigma$-modules satisfying \eqref{eq:polynomial weight}.
For $\la\in \cP_{M|N}$, let
\begin{equation}\label{eq:highest weight correspondence-2}
\La^\sigma_\la = \la'_1\de_{M+1} +\cdots +\la'_N\de_{M+N} + \mu_1\de_{1}+\cdots+\mu_M\de_{M},
\end{equation}
where $(\mu_1,\ldots,\mu_M) = (\max\{\la_1,N\}- N,\ldots,\max\{\la_M,N\}- N)$, and let $V^\sigma(\la)$ be the irreducible highest weight $\ov\U^\sigma$-module with highest weight $\La^\sigma_\la$. 
Then any irreducible $\ov\U^\sigma$-module in $\ov{\mc{O}}^\sigma_{\geq 0}$ is isomorphic to $V^\sigma(\la)$ for some $\la\in \cP_{M|N}$.

Let $\V^\sigma$ be the natural representation of $\ov\U^\sigma$, which is equal to $\W_1=\W_1(1)$ as a $\ov\U^\sigma$-module.
Then $(\mc{L}^\sigma,\mc{B}^\sigma):=(\mc{L}_1,\mc{B}_1)$  is also a crystal base of $\V^\sigma$, where $\mc{B}^\sigma$ is given by
\begin{equation*}
M+1\ \stackrel{^{M+1}}{\longrightarrow}\ \cdots  \ \stackrel{^{n-1}}{\longrightarrow} \ n \stackrel{^{0}}{\longrightarrow}\  1 \ \stackrel{^{1}}{\longrightarrow}\ \cdots  \ \stackrel{^{M-1}}{\longrightarrow} \ M.
\end{equation*}
Let $\I^\sigma=\I$ as a $\Z_2$-graded set with a linear ordering $M+1< \cdots < n < 1< \cdots < M$.
For a skew Young diagram  $\la/\mu$, let $SST^\sigma(\la/\mu)$ be the set of semistandard tableaux of shape $\la/\mu$ with letters in $\I^\sigma$. 

For $\la\in \cP_{M|N}$, there exists a connected $\ov I^\sigma$-colored oriented graph structure on $SST^\sigma(\la)$ with a unique highest weight element, which is a subgraph of $(\mc{B}^\sigma)^{\otimes m}$ with $m=\sum_{i\ge 1}\la_i$ (cf.~\eqref{admissible reading}). Then $V^\sigma(\la)$ has a crystal base $(L^\sigma(\la), B^\sigma(\la))$, where $(L^\sigma(\la), B^\sigma(\la))$ is defined as in \eqref{eq:(L,B)} replacing $\ov I$ with $\ov I^\sigma$, and the crystal $B^\sigma(\la)$ is isomorphic to $SST^\sigma(\la)$.

\begin{cor}\label{cor:restriction of KR} 
For $r, s\geq 1$ with $(s^r)\in \cP_{M|N}$, 
we have $\W_s^{(r)} \cong V^\sigma ((s^r))$ as a $\ov\U^\sigma$-module.
\end{cor}
\pf Since $\W_s^{(r)}$ is a finite-dimensional $\ov\U^\sigma$-module in $\ov{\mc{O}}^\sigma_{\geq 0}$, each irreducible component of $\W_s^{(r)}$ is $V^\sigma(\la)$ for some $\la\in \cP_{M|N}$. The character of $SST(\la)$ is the hook Schur polynomial $hs_\la(x_1,\ldots,x_M,y_1,\ldots,y_N)$ corresponding to $\la$, and it is also equal to the character of $SST^\sigma(\la)$ for $\la\in \cP_{M|N}$ (cf.~\cite{BR,BSS}).
This implies that the multiplicity of $V^\sigma(\la)$ in $\W_s^{(r)}$ is equal to the one of $V(\la)$ in $\W_s^{(r)}$ as a $\ov\U$-module for $\la\in \cP_{M|N}$. Hence $\W_s^{(r)} \cong V^\sigma ((s^r))$ as a $\ov\U^\sigma$-module by Proposition \ref{prop:W_s^{(r)}}.
\qed

\subsection{Description of $\te_0$ and $\tf_0$}
Let us describe the action of $\te_i$, $\tf_i$ for $i\in I$ on $B^{r,s}$.
We may identify $B^{r,s}=SST((s^r))$ as an $\ov I$-colored oriented graph by Theorem \ref{thm:crystal base of poly repn} and Theorem \ref{thm:crystal base of KR}. 
Hence it suffices to consider $\te_0$ and $\tf_0$.

For $\la\in \cP$, let $\la^\pi$ denote the skew Young diagram obtained from $\la$ by $180^\circ$-rotation.
Let $T$ be a semistandard tableau of shape $\la$ with letters in $\I_\varepsilon$ ($\varepsilon=0,1$).
We denote by $T^{\se}$ the unique tableau of shape $\la^\pi$ which is Knuth equivalent to $T$.
When $T$ is of shape $\la^\pi$, we denote by $T^{\nw}$ the unique tableau of shape $\la$ which is Knuth equivalent to $T$ (cf.~\cite{Ful}).

For $\la\in \cP_{M|N}$ and $T\in SST(\la)$, let $T_\varepsilon$ denote the subtableau of $T$ with letters in $\I_\varepsilon$ ($\varepsilon=0,1$). 
We may regard $T$ as the tableau obtained 
by putting together $T_0$ of shape $\mu$ and $T_1$ of shape $\la/\mu$ for some $\mu$, and write $T=T_0\ast T_1$. For $S\in SST^\sigma(\la)$, we also denote by $S_\varepsilon$ the subtableau of $S$ with letters in $\I_\varepsilon^\sigma$ ($\varepsilon=0,1$), and write $S=S_1\ast S_0$.

\begin{lem}\label{lem:sigma on SST}
For $r,s\geq 1$ with $(s^r)\in \cP_{M|N}$,
there exists a bijection
\begin{equation*}
\sigma : SST((s^r)) \longrightarrow SST^\sigma((s^r))
\end{equation*}
such that $\sigma(T)_0 = T_0^{\se}$ and $\sigma(T)_1 = T_1^{\nw}$, 
that is, $\sigma(T)  =T_1^{\nw}\ast T_0^{\se}$ for $T\in SST((s^r))$.
\end{lem}
\pf It is clear that $\sigma$ is well-defined.
Let $T\in SST((s^r))$ be given. 
Since $(T_0^{\se})^{\nw} = T_0$ and $(T_1^{\nw})^{\se} = T_1$, the map sending $S\in SST^\sigma(\la)$ to $S_0^{\nw}\ast S_1^{\se}\in SST((s^r))$ is the inverse of $\sigma$. Hence $\sigma$ is a bijection. \qed 

\begin{ex}\label{ex:switching}
{\rm
Suppose that $n=7$ and $\e=\e_{3|4}$, where $\I_0=\{1,2,3\}$ and $\I_1=\{4,5,6,7\}$.
Let
$$T=
\resizebox{.15\hsize}{!}
{\def\lr#1{\multicolumn{1}{|@{\hspace{.75ex}}c@{\hspace{.75ex}}|}{\raisebox{-.04ex}{$#1$}}}\raisebox{-.6ex}
{$\begin{array}{ccccc}
\cline{1-5}
\lr{\blue {\it 1}} & \lr{\blue {\it 1}} &  \lr{\blue {\it 1}} & \lr{\blue {\it 2}} & \lr{7} \\
\cline{1-5}
\lr{\blue {\it 2}} & \lr{\blue {\it 2}} & \lr{\blue {\it 3}} &  \lr{ 5} & \lr{ 7}  \\
\cline{1-5}
\lr{\blue {\it 3}} & \lr{ 4} & \lr{ 5} & \lr{ 6} & \lr{ 7} \\
\cline{1-5}
\end{array}$}}\quad \in SST((5^3)) = B^{3,5},$$\vskip 2mm
\noindent where the letters in $\I_0$ are denoted in blue (italic) for convenience. 
Then \vskip 2mm
$$T_0=
\resizebox{.15\hsize}{!}
{\def\lr#1{\multicolumn{1}{|@{\hspace{.75ex}}c@{\hspace{.75ex}}|}{\raisebox{-.04ex}{$#1$}}}\raisebox{-.6ex}
{$\begin{array}{ccccc}
\cline{1-4}
\lr{\blue {\it 1}} & \lr{\blue {\it 1}} &  \lr{\blue {\it 1}} & \lr{\blue {\it 2}} &   \\
\cline{1-4}
\lr{\blue {\it 2}} & \lr{\blue {\it 2}} & \lr{\blue {\it 3}} &      \\
\cline{1-3}
\lr{\blue {\it 3}} &   \\
\cline{1-1}
\end{array}$}}\quad\quad\quad\quad
T_0^{\se}=
\resizebox{.15\hsize}{!}
{\def\lr#1{\multicolumn{1}{|@{\hspace{.75ex}}c@{\hspace{.75ex}}|}{\raisebox{-.04ex}{$#1$}}}\raisebox{-.6ex}
{$\begin{array}{ccccc}
\cline{5-5}
 &   &   &  & \lr{\blue 1} \\
\cline{3-5}
 &   & \lr{\blue {\it 1}} &  \lr{\blue {\it 2}} & \lr{\blue {\it 2}}  \\
\cline{2-5}
& \lr{\blue {\it 1}} & \lr{\blue {\it 2}} & \lr{\blue {\it 3}} & \lr{\blue {\it 3}}   \\
\cline{2-5}
\end{array}$}}\quad ,$$\vskip 2mm
\noindent and  
$$T_1=
\resizebox{.15\hsize}{!}
{\def\lr#1{\multicolumn{1}{|@{\hspace{.75ex}}c@{\hspace{.75ex}}|}{\raisebox{-.04ex}{$#1$}}}\raisebox{-.6ex}
{$\begin{array}{ccccc}
\cline{5-5}
  &  &   &  & \lr{ 7}  \\
\cline{4-5}
  &  &   &  \lr{ 5} & \lr{ 7}  \\
\cline{2-5}
  & \lr{ 4} & \lr{ 5} & \lr{ 6} & \lr{ 7} \\
\cline{2-5}
\end{array}$}}\quad\quad\quad\quad
T_1^{\nw}=
\resizebox{.15\hsize}{!}
{\def\lr#1{\multicolumn{1}{|@{\hspace{.75ex}}c@{\hspace{.75ex}}|}{\raisebox{-.04ex}{$#1$}}}\raisebox{-.6ex}
{$\begin{array}{ccccc}
\cline{2-5}
  & \lr{ 4}  & \lr{ 5}   &  \lr{ 6} & \lr{ 7}  \\
\cline{2-5}
  &  \lr{ 5} & \lr{ 7} &  &  \\
\cline{2-3}
 & \lr{ 7} \\
 \cline{2-2}
\end{array}$}}\quad .$$
Hence we have \vskip 2mm
$$ \sigma(T) = 
\resizebox{.15\hsize}{!}
{\def\lr#1{\multicolumn{1}{|@{\hspace{.75ex}}c@{\hspace{.75ex}}|}{\raisebox{-.04ex}{$#1$}}}\raisebox{-.6ex}
{$\begin{array}{ccccc}
\cline{1-5}
\lr{ 4} & \lr{ 5} &  \lr{ 6} & \lr{ 7} & \lr{\blue {\it 1}} \\
\cline{1-5}
\lr{ 5} & \lr{ 7} & \lr{\blue {\it 1}} &  \lr{\blue {\it 2}} & \lr{\blue {\it 2}}  \\
\cline{1-5}
\lr{ 7} & \lr{\blue {\it 1}} & \lr{\blue {\it 2}} & \lr{\blue {\it 3}} & \lr{\blue {\it 3}} \\
\cline{1-5}
\end{array}$}}\quad .
$$

}
\end{ex}

\begin{rem}{\rm
The map $\sigma$ in Lemma \ref{lem:sigma on SST} can be described by exchanging the subtableaux  $T_0$ and $T_1$ applying {\em switching algorithm} in \cite{BSS}.
}
\end{rem}

Let us identify $B^{r,s}=SST((s^r))$. 
For $i\in \ov I^\sigma$, let $\te_i^\sigma$ and $\tf_i^\sigma$ denote the Kashiwara operators on the crystal $B^\sigma((s^r)) = SST^\sigma((s^r))$ of $\W_s^{(r)}$ as a $\ov\U^\sigma$-module.
Then we have the following combinatorial description of $\te_0$ and $\tf_0$ on $B^{r,s}$.

\begin{thm} \label{th:sigma}
Under the above hypothesis, we have
\begin{equation*}
\te_0  = \sigma^{-1} \circ \te_0^\sigma \circ \sigma,\quad
\tf_0  = \sigma^{-1} \circ \tf_0^\sigma \circ \sigma.
\end{equation*}
\end{thm}
\pf By Corollary \ref{cor:restriction of KR} and Theorem \ref{thm:uniqueness of crystal base}, there exists an isomorphism of $\ov I^\sigma$-colored oriented graphs
\begin{equation*}
\psi : B^{r,s} = SST((s^r)) \longrightarrow SST^\sigma ((s^r))=B^\sigma((s^r)).
\end{equation*}

We first claim that $\psi=\sigma$.
Let $T\in SST((s^r))$ be given. Assume that $T=T_0\ast T_1$ where $T_0$ is of shape $\mu$ and $T_1$ is of shape $(s^r)/\mu$ for some $\mu\in\cP$, and $\psi(T)=\psi(T)_1\ast\psi(T)_0$ where $\psi(T)_0$ is of shape $\nu^\pi$ and $\psi(T)_1$ is of shape $((s^r)/\nu)^\pi$ for some $\nu\in \cP$.

Since $\psi\circ \te_i =\te_i^\sigma \circ \psi$ and 
$\psi\circ \tf_i =\tf_i^\sigma \circ \psi$ for $i\in \ov I^\sigma$, 
the connected component of $T_0$ under $\te_i$ and $\tf_i$ for $i\in \{1,\ldots,M-1\}$ and the one of $\psi (T)_0$ under $\te^\sigma_i$ and $\tf^\sigma_i$ for $i\in \{1,\ldots,M-1\}$ are isomorphic as $\{1,\ldots,M-1\}$-colored oriented graphs, where the isomorphism maps $T$ to $\psi(T)_0$.
This implies that $\psi(T)_0$ is Knuth equivalent to $T_0$, and hence $\psi(T)_0=T_0^{\se}$ with $\mu=\nu$.
Applying similar arguments to $T_1$ with respect to 
$\te_i$ and $\tf_i$ for $i\in \{M+1,\ldots,n-1\}$, we have $\psi(T)_1 = T_1^{\nw}$. 
Therefore, $\psi(T)= T_1^{\nw}\ast T_0^{\se} = \sigma(T)$. 

Since $\psi=\sigma$ and $\psi$ is an isomorphism of $\ov I^\sigma$-colored oriented graphs, 
we have $\te_0  = \sigma^{-1} \circ \te_0^\sigma \circ \sigma$ and
$\tf_0  = \sigma^{-1} \circ \tf_0^\sigma \circ \sigma$. The proof completes. \qed

\begin{ex}\label{ex:f_0}
{\rm Let $T$ be as in Example \ref{ex:switching}.
Then   \vskip 2mm
$$ 
\sigma(T) = 
\resizebox{.15\hsize}{!}
{\def\lr#1{\multicolumn{1}{|@{\hspace{.75ex}}c@{\hspace{.75ex}}|}{\raisebox{-.04ex}{$#1$}}}\raisebox{-.6ex}
{$\begin{array}{ccccc}
\cline{1-5}
\lr{ 4} & \lr{ 5} &  \lr{ 6} & \lr{ 7} & \lr{\blue {\it 1}} \\
\cline{1-5}
\lr{ 5} & \lr{ 7} & \lr{\blue {\it 1}} &  \lr{\blue {\it 2}} & \lr{\blue {\it 2}}  \\
\cline{1-5}
\lr{ 7} & \lr{\blue {\it 1}} & \lr{\blue {\it 2}} & \lr{\blue {\it 3}} & \lr{\blue {\it 3}} \\
\cline{1-5}
\end{array}$}}\quad \quad \stackrel{0}{\longrightarrow}\quad\quad
\tf_0^\sigma\sigma(T) = 
\resizebox{.15\hsize}{!}
{\def\lr#1{\multicolumn{1}{|@{\hspace{.75ex}}c@{\hspace{.75ex}}|}{\raisebox{-.04ex}{$#1$}}}\raisebox{-.6ex}
{$\begin{array}{ccccc}
\cline{1-5}
\lr{ 4} & \lr{ 5} &  \lr{ 6} & \lr{ 7} & \lr{\blue {\it 1}} \\
\cline{1-5}
\lr{ 5} & \lr{ 7} & \lr{\blue {\it 1}} &  \lr{\blue {\it 2}} & \lr{\blue {\it 2}}  \\
\cline{1-5}
\lr{\blue {\it 1}} & \lr{\blue {\it 1}} & \lr{\blue {\it 2}} & \lr{\blue {\it 3}} & \lr{\blue {\it 3}} \\
\cline{1-5}
\end{array}$}}\quad 
$$\vskip 2mm
\noindent by applying \eqref{eq:tensor product rule for odd -} to $\iota(\sigma(T))$ \eqref{admissible reading}, for example, $\iota(\sigma(T))=(123)(723)(613)(572)(457)$ given by reading the letters of $\sigma(T)$ column by column from right to left and from top to bottom in each column.
On the other hand,
$$\sigma^{-1}(\tf_0^\sigma\sigma(T))=
\resizebox{.15\hsize}{!}
{\def\lr#1{\multicolumn{1}{|@{\hspace{.75ex}}c@{\hspace{.75ex}}|}{\raisebox{-.04ex}{$#1$}}}\raisebox{-.6ex}
{$\begin{array}{ccccc}
\cline{1-5}
\lr{\blue {\it 1}} & \lr{\blue {\it 1}} &  \lr{\blue {\it 1}} & \lr{\blue {\it 1}} & \lr{\blue {\it 2}} \\
\cline{1-5}
\lr{\blue {\it 2}} & \lr{\blue {\it 2}} & \lr{\blue {\it 3}} &  \lr{ 5} & \lr{ 7}  \\
\cline{1-5}
\lr{\blue {\it 3}} & \lr{ 4} & \lr{ 5} & \lr{ 6} & \lr{ 7} \\
\cline{1-5}
\end{array}$}}\quad .$$\vskip 2mm
\noindent Therefore, \vskip 2mm
$$ 
\resizebox{.15\hsize}{!}
{\def\lr#1{\multicolumn{1}{|@{\hspace{.75ex}}c@{\hspace{.75ex}}|}{\raisebox{-.04ex}{$#1$}}}\raisebox{-.6ex}
{$\begin{array}{ccccc}
\cline{1-5}
\lr{\blue {\it 1}} & \lr{\blue {\it 1}} &  \lr{\blue {\it 1}} & \lr{\blue {\it 2}} & \lr{7} \\
\cline{1-5}
\lr{\blue {\it 2}} & \lr{\blue {\it 2}} & \lr{\blue {\it 3}} &  \lr{ 5} & \lr{ 7}  \\
\cline{1-5}
\lr{\blue {\it 3}} & \lr{ 4} & \lr{ 5} & \lr{ 6} & \lr{ 7} \\
\cline{1-5}
\end{array}$}}\quad\quad  \stackrel{0}{\longrightarrow} \quad\quad
\resizebox{.15\hsize}{!}
{\def\lr#1{\multicolumn{1}{|@{\hspace{.75ex}}c@{\hspace{.75ex}}|}{\raisebox{-.04ex}{$#1$}}}\raisebox{-.6ex}
{$\begin{array}{ccccc}
\cline{1-5}
\lr{\blue {\it 1}} & \lr{\blue {\it 1}} &  \lr{\blue {\it 1}} & \lr{\blue {\it 1}} & \lr{\blue {\it 2}} \\
\cline{1-5}
\lr{\blue {\it 2}} & \lr{\blue {\it 2}} & \lr{\blue {\it 3}} &  \lr{ 5} & \lr{ 7} \\
\cline{1-5}
\lr{\blue {\it 3}} & \lr{ 4} & \lr{ 5} & \lr{ 6} & \lr{ 7} \\
\cline{1-5}
\end{array}$}}\quad .$$\vskip 3mm

We also include two crystal graphs in Figures 1 and 2.

\begin{figure}[h]
\includegraphics[scale=.75]{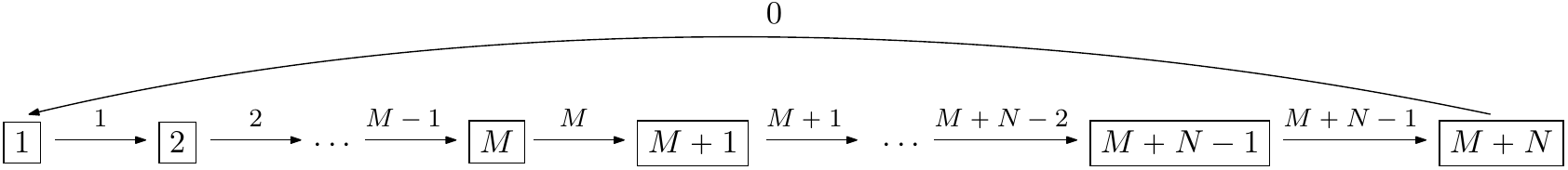}
\caption{Crystal graph of $B^{1,1}$ for $\U(\epsilon_{M|N})$}
\end{figure}

\begin{figure}[h]
\includegraphics[scale=.75]{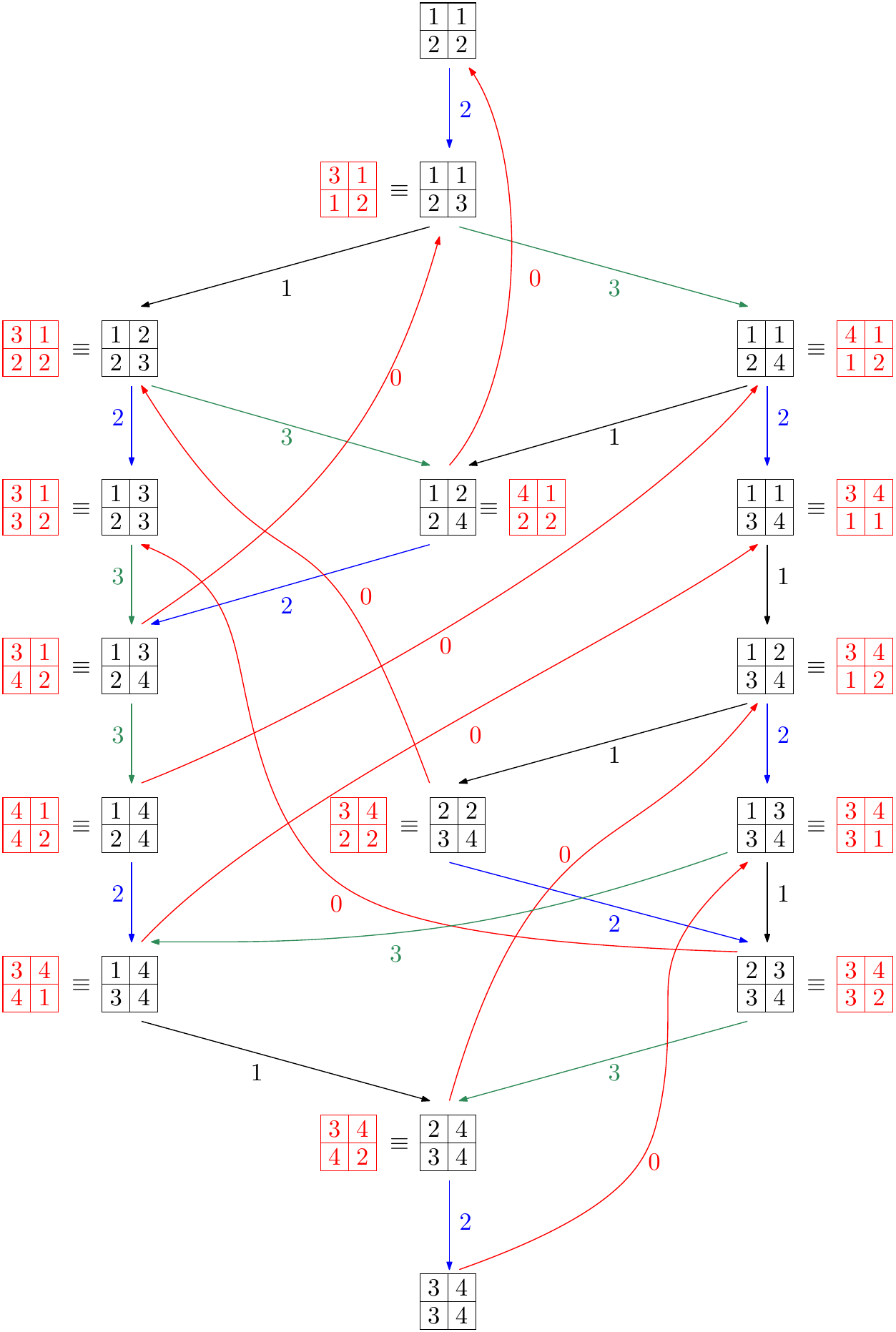}
\caption{Crystal graph of $B^{2,2}$ for $\U(\epsilon_{2|2})$ (the red tableaux denote the images under $\sigma$)}
\end{figure}
}
\end{ex}

\section{Combinatorial $R$ matrix}\label{sec:combR}

\subsection{Genuine highest weight vectors in $B^{r_1,s_1}\ot B^{r_2,s_2}$}
Let $(s_1^{r_1}), (s_2^{r_2})\in \cP_{M|N}$ be given.
We set $r=\min(r_1,r_2)$ and $(a)_+=\max(a,0)$ for $a\in\Z$.
For a partition $\nu\subset(s_1^{\min(r_1,M)})$, $\overline{\nu}$ stands for the one 
adding $s_2-s_1$ columns of height $\min(r_1,M)$ to the left of $\nu$ if $s_1\le s_2$, and 
removing leftmost $s_1-s_2$ columns from $\nu$ if $s_1\ge s_2$.

Let $T_1\ot T_2 \in B^{r_1,s_1}\ot B^{r_2,s_2}$ be given.
Let $col(T_2)$ denote the word $w_1\cdots w_{l}$ given by reading the letters in $T_2$ column by column from right to left and top to bottom in each column. We let $col(T_2)\rightarrow T_1$ be the tableau obtained by applying column insertion of $col(T_2)$ to $T_1$ (starting from $w_1$) (cf.~\cite{Ful,KK}). 
Recall that $col(T_2)\rightarrow T_1$ is a semistandard tableau whose $\ov{I}$-component is isomorphic to that of $T_1\otimes T_2$ (see \cite[Section 4.5]{BKK}). Here an $\ov{I}$-component in a crystal means a connected component under $\te_i$, $\tf_i$ for $i\in \ov{I}$. We also say that $T_1\ot T_2$ is a genuine highest weight vector if $col(T_2)\rightarrow T_1$ is so.

Now, let us parametrize genuine highest weight vectors in $B^{r_1,s_1}\ot B^{r_2,s_2}$. 
By Corollary \ref{cor:LR rule for poly} and Proposition \ref{prop:W_s^{(r)}}, 
the crystal $B^{r_1,s_1}\ot B^{r_2,s_2}$ 
has a multiplicity-free decomposition into $\ov{I}$-components parametrized by certain partitions in $\cP_{M|N}$.
Such a partition, which we denote by $\widehat{\la}$, is constructed from a partition $\la$ with $\ell(\la)\le r$ by placing $\la$ right to $(s_2^{r_2})$
and then $\la^c$ below it, where $\la^c$ is the partition obtained by rotating
$(s_1^{r_1})/\la$ by $180^\circ$. 

\begin{prop} \label{pr:genuine highest} 
The partition $\la$ given above and the corresponding genuine highest weight vector $T_1\ot T_2 \in B^{r_1,s_1}\ot B^{r_2,s_2}$ with weight $\widehat{\la}$ are given by one of the followings:
\begin{enumerate}
\item When $r_1,r_2>M$, fix partitions $\mu,\nu$ such that 
$$(((s_1-s_2)_+)^{r-M})\subset
\mu\subset(s_1^{r-M}),\quad
\nu\subset(s_1^M),\quad \mu_1\le \min(\nu_M,N-s_2).$$ 
{\rm (i)} $\la$ is constructed from $\mu,\nu$ by adding $\mu$ below $\nu$. 
{\rm (ii)} $T_1$ is divided into three parts: In the first $M$ rows, letters are $i$ in the $i$-th row. 
Below the $M$-th row, in the lower right part with the shape obtained by rotating 
$\mu$ by $180^\circ$, letters range as $M+s_2+1,M+s_2+2,\ldots$ horizontally, and in 
the other part, letters range as $M+1,M+2,\ldots$ horizontally. 
{\rm (iii)} $T_2$ is divided into two parts: In the upper left part corresponding to $\overline{\nu}$, letters are 
$i$ in the $i$-th row. In the other part, letters range as $M+1,M+2,\ldots$ horizontally 
in each row.
\item When $r_1,r_2\le M$, fix a partition $\nu$  such that 
$\nu\subset(s_1^{r})$. 
{\rm (i)} $\la$ is given by $\nu$.
{\rm (ii)} $T_1$ is the genuine highest weight vector. 
{\rm (iii)} $T_2$ is divided into three parts $\overline{\nu}\subset\eta\subset(s_2^{r_2})$: In the upper left part corresponding to $\overline{\nu}$, letters are $i$ in the $i$-th row.
In the middle part $\eta/\overline{\nu}$, letters range as $r_1+1,r_1+2,\ldots$ up to
at most $M$ vertically. Finally, in the lower right part $(s_2^{r_2})/\eta$, 
letters range as $M+1,M+2,\ldots$ horizontally.
\item When $r_1>M,r_2\le M$ or $r_1\le M,r_2>M$, fix a partition $\nu$ as in {\rm (2)}.
{\rm (i)} $\la$ is given by $\nu$.
{\rm (ii)} $T_1$ is the genuine highest weight vector. 
{\rm (iii)} $T_2$ is the same as in {\rm (1)} when $r_1>M,r_2\le M$, and as in {\rm (2)} when $r_1\le M,r_2>M$.
\end{enumerate}
Here we assume that $\mu$ and $\nu$ are chosen only when $\widehat{\la}\in \cP_{M|N}$.
\end{prop}

\pf
It suffices to check that the tableau $col(T_2)\rightarrow T_1$ is a genuine 
highest weight vector with highest weight $\widehat{\la}$.
\qed

The goal of this subsection is to show that $B^{r_1,s_1}\ot B^{r_2,s_2}$ is connected. For this, we need to consider the column insertion of $T_2$ into $T_1$ in more details to see how the operators $\te_0$ and $\tf_0$ act on each $\ov{I}$-component $T_1\ot T_2$ in $B^{r_1,s_1}\ot B^{r_2,s_2}$. 

Suppose that $T_1\ot T_2$ is a genuine highest weight vector with highest weight $\widehat{\la}$
as described in Proposition \ref{pr:genuine highest}. We also keep the other notations in Proposition \ref{pr:genuine highest}.
Let $\ov{\mu}$ be the partition given by
adding $s_2-s_1$ columns of height $(r-M)_+$ to the left of $\mu$ if $s_1\le s_2$, and 
removing leftmost $s_1-s_2$ columns from $\mu$ if $s_1\ge s_2$.
We assume that $\mu$ is empty unless $r_1, r_2>M$.
Let $\td{\la}$ be the partition given by placing $\ov{\mu}$ below $\ov{\nu}$. 
Note that $\td{\la}$ is also obtained from $\la$ by adding
$s_2-s_1$ (resp. removing $s_1-s_2$) columns of height $r$ if $s_1 \le s_2$
(resp. $s_1 \ge s_2$), and $\widehat{\la}$ is equal to the one given by placing $\td{\la}$ to the right of $(s_1^{r_1})$ and the compliment of $\td{\la}$ in $(s_2^{r_2})$ below it.

We define $T_2^\circ$ to be a tableau of shape $(s_2^{r_2})$ with letters in $\N$, where the letter at each node in $(s_2^{r_2})$ is determined by the column index $a$ (enumerated from the left) if a node at the $a$-th column is increased after the insertion of the letter of $T_2$ at the same position. Note that $T_2^\circ$ is not semistandard, but the anti-diagonal flip of $T_2^{\circ}$ is semistandard (see \cite[Section 5]{KK}).

\begin{lem}\label{lem:recording for genuine}
Under the above hypothesis,  
$T_2^{\circ}$ is determined by the rule, where
letters range $\cdots, s_1+2, s_1+1$ horizontally from right to left in $\td{\la}$, and letters range $\cdots, 2, 1$ horizontally from right to left elsewhere. 
\end{lem}
\pf
It follows directly from the description of $T_2$ in Proposition \ref{pr:genuine highest} and column insertion.
\qed
 
Let $\zeta$ be the partition given by removing a node from the rightmost column of $\la$.
Suppose that $\widehat{\zeta}\in \cP_{M|N}$. 
We first show that the $\ov{I}$-component in $B^{r_1,s_1}\ot B^{r_2,s_2}$ with highest weight $\widehat{\la}$ is connected to that of $\widehat{\zeta}$ when $\max(r_1,r_2)>M$ or  $r_1+r_2\leq M$ 
in the following two lemmas.

First suppose that either $r_1>M$ or $r_2>M$.
Let $i_1,i_2,\ldots,i_s$ ($s=N+\nu'_{\nu_1}+\la_1-s_1-1$) be the sequence given by
\begin{equation}\label{eq:seq of indices-1}
0,M+N-1,M+N-2,\ldots,M+s_1-\la_1+1,1,2,\ldots,\nu'_{\nu_1}-1.
\end{equation} 
Set $T'_1\ot T'_2=\te_{i_s}\ldots\te_{i_2}\te_{i_1}(T_1\ot T_2)$. 
Let us say that $T_1\ot T_2$ belongs to case (E) if and only if $r_1, r_2>M$ and
$\nu$ is of rectangular shape with $\nu_M=\mu_1$.

\begin{lem} \label{lem:genuine highest-1}
Under the above hypothesis, we have the following. 
\begin{itemize}
\item[(1)] If it is not case (E), then $T'_1\ot T'_2$ is a genuine highest weight vector with weight $\widehat{\zeta}$. 

\item[(2)] If it is case (E), then $T'_1\ot T'_2$ is almost the same as a genuine highest 
weight vector with weight $\widehat{\zeta}$ except that the letter at the rightmost 
node in the $M$-th row of $col(T'_2)\rightarrow T'_1$ is $M+s_2+\mu_1$.
\end{itemize}
In either case, the application of $\te_j$ for $j=i_1,i_2,\ldots,i_s$ always takes
place on the second component, that is, $T'_1=T_1$ and $T'_2=\te_{i_s}\ldots\te_{i_2}\te_{i_1}T_2$.
\end{lem}

\pf Suppose that $r_1, r_2>M$. 
Let $x$ be the node in $(s_2^{r_2})$ corresponding to the rightmost corner of $\widehat{\la}$ and let $y$ be the node in $(s_2^{r_2})$ corresponding to the rightmost corner in $\ov{\mu}$ (regarded as a subdiagram of $\td{\la}$ placed below $\ov{\nu}$). Let $z$ be a position of any node in $(s_2^{r_2})$ below $x$ in the same column.

\begin{figure}[h]
\includegraphics[scale=.75]{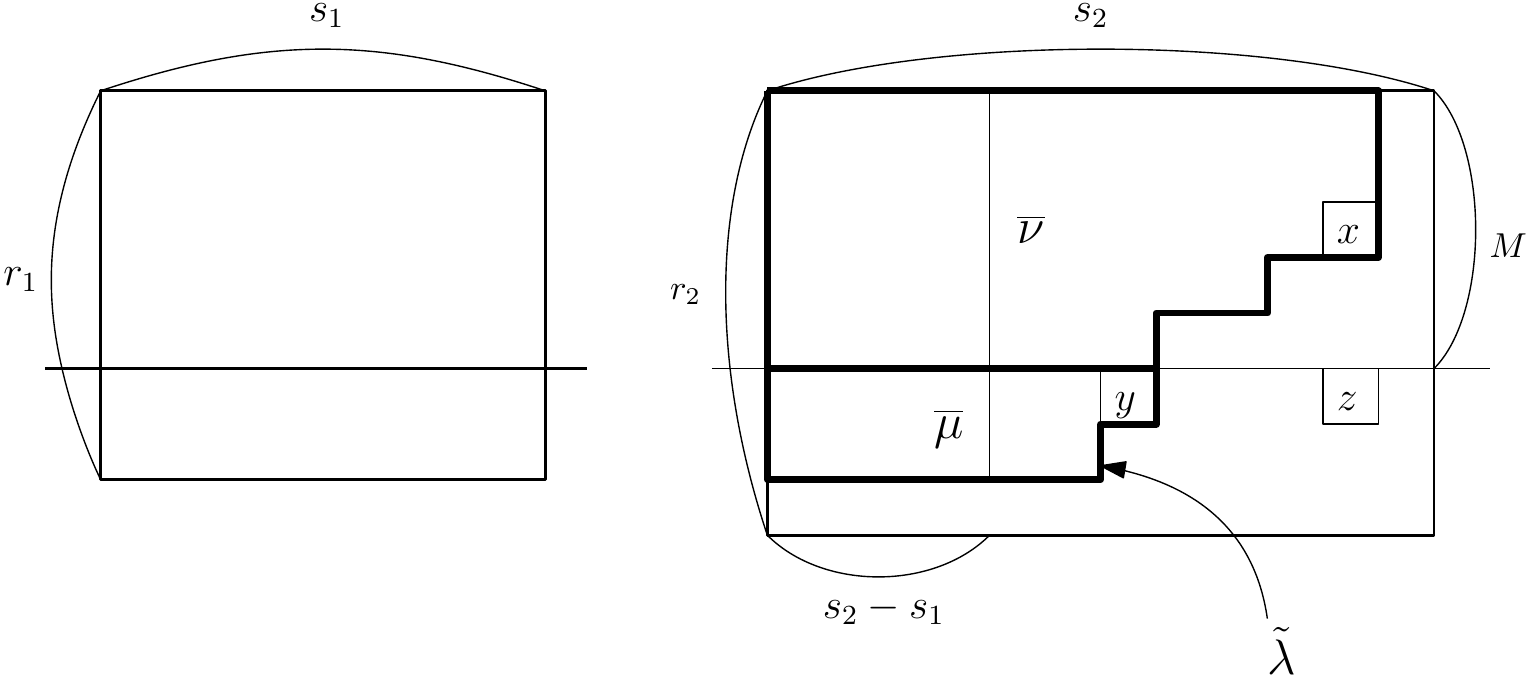}
\caption{Partitions $\ov{\mu}$, $\ov{\nu}$ and $\tilde{\lambda}$ in $(s_2^{r_2})$}
\end{figure}

{\em Case 1.} Suppose that $T_1\ot T_2$ does not belong to case (E).
Then $x\neq y$, and the letters at $x$ and $y$ in $T_2^\circ$ are 
$s_1+1$ by Lemma \ref{lem:recording for genuine}. 
Since $y$ is located strictly to the southwest of $x$, the letter at $z$ is strictly smaller than $s_1+1$.

Let us consider $T'_1\ot T'_2=\te_{i_s}\ldots\te_{i_2}\te_{i_1}(T_1\ot T_2)$. 
By using the description of $T_1$ and $T_2$ in Proposition \ref{pr:genuine highest}(1) and Theorem \ref{th:sigma}, 
we see that $\te_j$ for $j=i_1,i_2,\ldots,i_s$ acts only on the second component in the tensor product, that is,
$T'_1=T_1$ and $T'_2=\te_{i_s}\ldots\te_{i_2}\te_{i_1} T_2$. Furthermore, we observe that $T'_2$ is equal to $T_2$ except that the letters starting at $x$ are given by $M+1, M+2,\cdots$ horizontally.

By Lemma \ref{lem:recording for genuine} and applying column insertion of $T'_2$ into $T_1$, it is not difficult to see that $(T'_2)^\circ$, which is defined in the same way as $T_2^\circ$ with respect to $T'_2$, is equal to $T_2^\circ$ 
except that the letters starting at $x$ are given by
$\cdots, s_1+2, s_1+1, u$ horizontally, where $u=s_1-\la_1+1$. 
This implies that $col(T'_2)\rightarrow T_1$ is equal to $H_{\widehat{\zeta}}$, the genuine highest weight vector with weight $\widehat{\zeta}$.

{\em Case 2}. Suppose that $T_1\ot T_2$ belongs to case (E). 
In this case, we have $x=y$, where in $T_2^\circ$ the letter at $x$ is $s_1+1$, and the letters below and to the right of $x$ are no more than $s_1$.

Let $x'$ be the node in $(s_2^{r_2})$ corresponding to the rightmost corner of $\ov{\nu}$.
As in {\em Case 1}, we have $T'_1=T_1$ and $T'_2=\te_{i_s}\ldots\te_{i_2}\te_{i_1} T_2$, where  $T'_2$ is equal to $T_2$ except that the letters starting at $x'$ are given by $M+1, M+2,\cdots$ horizontally.
Moreover, $(T'_2)^\circ$ is equal to $T_2^\circ$ except that 
the letters starting at $x$ are given by
$\cdots, s_1+2, s_1+1, u$ horizontally.
But then $col(T'_2)\rightarrow T_1$ is equal to a tableau obtained from 
$H_{\widehat{\zeta}}$ by replacing the letter at $x'$ the rightmost node in the $M$-th row with $M+s_2+\mu_1$.

The proof for the case when either $r_1>M, r_2\leq M$ or $r_1\leq M, r_2>M$ are similar, and we leave it to the reader.
\qed

Next, suppose that $r_1+r_2\leq M$. In particular we have $r<M$. 
Let $j_1,j_2,\ldots,j_t$ ($t=M+N+2\la'_{\la_1}-r_1-r_2-1$) be the sequence given by
\begin{equation}\label{eq:seq of indices-2}
0,M+N-1,\ldots, M+1,M,\ldots, r_1+r_2-\la'_{\la_1}+1,1,2,\ldots,\la'_{\la_1}-1.
\end{equation} 
Set $T''_1\ot T''_2=\te_{j_t}\ldots\te_{j_2}\te_{j_1}(T_1\ot T_2)$. Then we can check the following in a straightforward manner.

\begin{lem} \label{lem:genuine highest-2}
Under the above hypothesis, $T''_1\ot T''_2$ is a genuine highest weight vector with weight $\widehat{\zeta}$. The application of $\te_j$ for $j=j_1,j_2,\ldots,j_t$ always takes
place on the second component, that is, $T''_1=T_1$ and $T''_2=\te_{j_t}\ldots\te_{j_2}\te_{j_1}T_2$.
\end{lem}

\begin{ex}{\rm
(1) Let $M=2,N=5$ and

\begin{equation*}
T_1\ot T_2 = \vcenter{
\tableau[sby]{1&1&1&1\\2&2&2&2\\3&4&5&6\\3&5&6&7}
}
\otimes
\vcenter{
\tableau[sby]{1&1\\2&3\\3&4\\3&4}
}\,.
\end{equation*}
We have
\begin{equation*}
(col(T_2)\rightarrow T_1)=
\vcenter{
\tableau[sby]{1&1&1&1&1&1\\2&2&2&2&2\\3&4&5&6&7\\3&4&5&6\\3&4\\3\\3}
}\quad\quad\quad
T_2^\circ = \vcenter{
\tableau[sby]{6&5\\5&1\\5&1\\2&1}
}\,,
\end{equation*}
\noindent where $\td{\la}=(2,1,1)$ with $\ov{\nu}=(2,1)$ and $\ov{\mu}=(1)$
following the notations in the proof of Lemma \ref{lem:genuine highest-1}. 
It is not case (E). Then applying a sequence of $\te_i$'s according to \eqref{eq:seq of indices-1}
\begin{equation*}
\begin{split}
& T_2 =
\vcenter{
\tableau[sby]{1&1\\2&3\\3&4\\3&4}
}
\equiv 
\vcenter{
\tableau[sby]{3&4\\3&4\\3&1\\1&2}
}\quad \\ \\ & \stackrel{0}{\longleftarrow} \quad
\vcenter{
\tableau[sby]{3&4\\3&4\\3&1\\7&2}
}\quad  \stackrel{6}{\longleftarrow} \quad
\vcenter{
\tableau[sby]{3&4\\3&4\\3&1\\6&2}
}\quad  \stackrel{5}{\longleftarrow} \quad
\vcenter{
\tableau[sby]{3&4\\3&4\\3&1\\5&2}
}\quad\stackrel{4}{\longleftarrow} \quad
\vcenter{
\tableau[sby]{3&4\\3&4\\3&1\\4&2}
}\quad  \stackrel{3}{\longleftarrow} \quad
\vcenter{
\tableau[sby]{3&4\\3&4\\3&1\\3&2}
} \equiv
\vcenter{
\tableau[sby]{1&3\\2&3\\3&4\\3&4}
} = T'_2,
\end{split}
\end{equation*}
and
\begin{equation*}
(col(T'_2)\rightarrow T_1)=
\vcenter{
\tableau[sby]{1&1&1&1&1\\2&2&2&2&2\\3&4&5&6&7\\3&4&5&6\\3&4\\3\\3\\3}
}
\quad\quad\quad
(T'_2)^\circ = \vcenter{
\tableau[sby]{5&1\\5&1\\5&1\\2&1}
}\,.
\end{equation*}

(2) Let $T_1\ot T'_2$ be as in (1). In this case, we have $\ov{\nu}=(1,1)$, $\ov{\mu}=(1)$, and $\td{\la}=(1,1,1)$. It belongs to case (E). Then 
\begin{equation*}
\begin{split}
& T'_2  = 
\vcenter{
\tableau[sby]{1&3\\2&3\\3&4\\3&4}
}
\equiv 
\vcenter{
\tableau[sby]{3&4\\3&4\\3&1\\3&2}
}\quad \\ \\ & \stackrel{0}{\longleftarrow} \quad
\vcenter{
\tableau[sby]{3&4\\3&4\\3&7\\3&2}
} \quad \stackrel{6}{\longleftarrow} \quad
\vcenter{
\tableau[sby]{3&4\\3&4\\3&6\\3&2}
}\quad  
\stackrel{5}{\longleftarrow} \quad
\vcenter{
\tableau[sby]{3&4\\3&4\\3&5\\3&2}
}\quad \stackrel{4}{\longleftarrow} \quad
\vcenter{
\tableau[sby]{3&4\\3&4\\3&4\\3&2}
}\quad \stackrel{1}{\longleftarrow} \quad
\vcenter{
\tableau[sby]{3&4\\3&4\\3&4\\4&1}
} \equiv
\vcenter{
\tableau[sby]{1&3\\3&4\\3&4\\3&4}
} = T''_2,
\end{split}
\end{equation*}
and
\begin{equation*}
(col(T''_2)\rightarrow T_1)=
\vcenter{
\tableau[sby]{1&1&1&1&1\\2&2&2&2&7\\3&4&5&6\\3&4&5&6\\3&4\\3&4\\3\\3}
}
\quad\quad\quad
(T''_2)^\circ = \vcenter{
\tableau[sby]{5&1\\5&1\\2&1\\2&1}
}\,.
\end{equation*}

(3) Let $M=6$, $N=2$, and
\begin{equation*}
T_1\ot T_2 = \vcenter{
\tableau[sby]{1&1&1\\2&2&2\\3&3&3}
}
\ot
\vcenter{
\tableau[sby]{1&1&4\\2&2&5\\3&4&6}
}\,.
\end{equation*}
We have 
\begin{equation*}
(col(T_2)\rightarrow T_1)=
\vcenter{
\tableau[sby]{1&1&1&1&1\\2&2&2&2&2\\3&3&3&3\\4&4\\5\\6}
}\quad\quad\quad
T_2^\circ = \vcenter{
\tableau[sby]{5&4&1\\5&4&1\\4&2&1}
}\,,
\end{equation*}
where $\td{\la}=\ov{\nu}=(2,2,1)$. Then applying a sequence of $\te_i$'s according to \eqref{eq:seq of indices-2}
\begin{equation*}
\begin{split}
& T_2= \vcenter{
\tableau[sby]{1&1&4\\2&2&5\\3&4&6}
}
\quad \stackrel{0}{\longleftarrow} \quad
\vcenter{
\tableau[sby]{8&1&4\\2&2&5\\3&4&6}
}\quad \equiv \quad
\vcenter{
\tableau[sby]{1&2&4\\2&4&5\\3&6&8}
}\\ \\
& \quad \stackrel{7}{\longleftarrow} \quad
\vcenter{
\tableau[sby]{1&2&4\\2&4&5\\3&6&7}
}
\quad \stackrel{6}{\longleftarrow} \quad
\vcenter{
\tableau[sby]{1&2&4\\2&4&5\\3&6&6}
}
\quad \stackrel{5}{\longleftarrow} \quad
\vcenter{
\tableau[sby]{1&2&4\\2&4&5\\3&5&6}
}
\quad \stackrel{1}{\longleftarrow} \quad
\vcenter{
\tableau[sby]{1&1&4\\2&4&5\\3&5&6}
} = T''_2,
\end{split}
\end{equation*}
and
\begin{equation*}
(col(T''_2)\rightarrow T_1) = 
\vcenter{
\tableau[sby]{1&1&1&1&1\\2&2&2&2\\3&3&3&3\\4&4\\5&5\\6}
}\,.
\end{equation*}

}
\end{ex}

Now we can prove the following.

\begin{thm}\label{th:connectedness of tensor product}
The crystal $B^{r_1,s_1}\ot B^{r_2,s_2}$ is connected.
\end{thm}
\pf 
{\em Case 1}. Suppose that $r_1>M$ or $r_2>M$. 
Let $\la_{min}$ be as in Proposition \ref{pr:genuine highest} 
such that $\widehat{\la_{min}}\in \cP_{M|N}$ and the number of nodes is minimal.
In this case, we have  $\la_{min}=(((s_1-s_2)_+)^r)$.

Let $T_1\ot T_2\in B^{r_1,s_1}\ot B^{r_2,s_2}$ be a genuine highest weight vector with weight $\widehat{\la}$ with $\la\neq \la_{min}$.
Let $\zeta$ be the partition given by removing a node from the rightmost column of $\la$. 
Then we have $\widehat{\zeta}\in \cP_{M|N}$. 
By Lemma \ref{lem:genuine highest-1}, 
$T_1\ot T_2$ is connected to the genuine highest weight vector with weight $\widehat{\zeta}$ in $B^{r_1,s_1}\ot B^{r_2,s_2}$. 
We may repeat this step to see that $T_1\ot T_2$ is connected to the genuine highest weight vector with weight $\widehat{\la_{min}}$. Hence $B^{r_1,s_1}\ot B^{r_2,s_2}$ is connected.

{\em Case 2}. Suppose that $r_1, r_2\le  M$ and $s_1, s_2\le N$ with $r_1+r_2> M$.
The proof is the same as that of {\em Case 1}.

{\em Case 3}. Suppose that $r_1, r_2\le M$ and $s_1, s_2\le N$ with $r_1+r_2\leq M$.
We apply the same argument as in {\em Case 1} now using Lemma \ref{lem:genuine highest-2} instead of Lemma \ref{lem:genuine highest-1}. Then we conclude that $B^{r_1,s_1}\ot B^{r_2,s_2}$ is connected.

{\em Case 4}. Suppose that $r_1, r_2\le M$ and $s_i>N$ for some $i=1,2$. Consider $\sigma(T_1)\ot \sigma(T_2)$ as an element of a crystal over 
$\U (\underbrace{1,\cdots,1}_{N},\underbrace{0\cdots,0}_{M})$. Then we may apply the result in {\em Case 1} to conclude that $B^{r_1,s_1}\ot B^{r_2,s_2}$ is connected. 
\qed

\subsection{Combinatorial $R$ matrix and energy function} \label{subsec:comb R}
Let $(s_1^{r_1}), (s_2^{r_2})\in \cP_{M|N}$ be given.
For $B^{r_1,s_1}\ot B^{r_2,s_2}$ we call a bijection 
$R: B^{r_1,s_1}\ot B^{r_2,s_2} \to B^{r_2,s_2}\ot B^{r_1,s_1}$ {\em combinatorial $R$ matrix},
if it commutes with $\te_i,\tf_i$ for any $i\in I$. Since $B^{r_1,s_1}\ot B^{r_2,s_2}$
is connected by Theorem \ref{th:connectedness of tensor product}, it is unique if it exists.
An integer valued function $H$ on $B^{r_1,s_1}\ot B^{r_2,s_2}$ is called {\em energy function}
if $H$ is constant on $\ov{I}$-components, and for
$b=T_1\otimes T_2\in B^{r_1,s_1}\ot B^{r_2,s_2}$,
\begin{align} \label{def H}
  H(\te_0b) &= H(b) +
  \begin{cases}
   1 & \text{in case LL}, \\
   0 &\text{in case LR or RL}, \\
   -1 &\text{in case RR},
  \end{cases}
\end{align}
where in case LL, when $\te_0$ is applied to $T_1\otimes T_2$ and to
$R(T_1\otimes T_2)=\tilde{T}_2\otimes \tilde{T}_1$, it acts on the left factor both times, 
in case RR $\te_0$ acts on the right factor both times, etc.
The existence of the combinatorial $R$ matrix and the energy function in our case can be 
shown in a similar way to \cite[Proposition 4.3.1]{KMN} or \cite[Proposition 2.6]{O:Memoirs} 
from the existence of an $R$ matrix for $\W^{(r_1)}_{s_1}\ot\W^{(r_2)}_{s_2}$ 
by fusion construction.

\begin{lem} \label{lem:comb R}
Let $T_1\ot T_2$ be a genuine highest weight vector with weight $\widehat{\la}$ for $\la$ as in Proposition \ref{pr:genuine highest}.
Suppose that $T_1\otimes T_2$ is sent to $\td{T}_2\otimes\td{T}_1$ 
by the combinatorial $R$ matrix.
Then $\td{T}_2\ot\td{T}_1$ is the genuine highest weight vector with weight $\widehat{\tilde{\la}}$.
\end{lem}
\pf It follows from the fact $\widehat{\tilde{\la}}=\widehat{\la}$.
\qed

\begin{ex}{\rm
\noindent(1) Let $M=2,N\ge5$. By the combinatorial $R$ matrix,\vskip 2mm
\begin{equation*}
\hspace{5mm}\vcenter{
\tableau[sby]{1&1&1&1\\2&2&2&2\\3&4&5&6\\3&5&6&7}
}
\otimes
\vcenter{
\tableau[sby]{1&1\\2&3\\3&4\\3&4}
}
\quad\text{is sent to}\quad
\vcenter{
\tableau[sby]{1&1\\2&2\\3&4\\3&7}
}
\otimes
\vcenter{
\tableau[sby]{1&1&1&1\\2&2&2&3\\3&4&5&6\\3&4&5&6}
}\,.
\end{equation*}\vskip 2mm
Both hand sides belong to the case (1) in Proposition \ref{pr:genuine highest}.
The corresponding partitions are $\la=(4,3,3,2),\tilde{\la}=(2,1,1)$.

\noindent(2) Let $M=4,N\ge3$. By the combinatorial $R$ matrix, \vskip 2mm
\begin{equation*}
\hspace{5mm}\vcenter{
\tableau[sby]{1&1&1\\2&2&2\\3&3&3}
}
\otimes
\vcenter{
\tableau[sby]{1&1&1&1\\2&2&4&4\\3&4&5&6\\4&5&6&7}
}
\quad\text{is sent to}\quad
\vcenter{
\tableau[sby]{1&1&1&1\\2&2&2&2\\3&3&3&3\\4&4&4&4}
}
\otimes
\vcenter{
\tableau[sby]{1&1&1\\2&5&6\\5&6&7}
}\,.
\end{equation*}\vskip 2mm
Both hand sides belong to the case (2) in Proposition \ref{pr:genuine highest}.
The corresponding partitions are $\la=(3,1),\tilde{\la}=(4,2,1)$.

\noindent(3) Let $M=3,N\ge3$. By the combinatorial $R$ matrix,\vskip 2mm
\begin{equation*}
\hspace{5mm}\vcenter{
\tableau[sby]{1&1&1&1\\2&2&2&2\\3&3&3&3\\4&5&6&7}
}
\otimes
\vcenter{
\tableau[sby]{1&1&4\\4&5&6}
}
\quad\text{is sent to}\quad
\vcenter{
\tableau[sby]{1&1&1\\2&2&2}
}
\otimes
\vcenter{
\tableau[sby]{1&1&1&3\\2&3&3&4\\3&4&5&6\\4&5&6&7}
}\,.
\end{equation*}\vskip 2mm
Both hand sides belong to the case (3) in Proposition \ref{pr:genuine highest}.
The corresponding partitions are $\la=(3,1),\tilde{\la}=(2)$.}
\end{ex}

Lemma \ref{lem:comb R} shows the image of the combinatorial $R$ matrix on the 
genuine highest weight elements. This can be generalized to all elements in 
$B^{r_1,s_1}\ot B^{r_2,s_2}$ in terms of insertion algorithms.

\begin{thm}\label{th:comb R}
The combinatorial $R$ matrix
\[
R:B^{r_1,s_1}\ot B^{r_2,s_2}\longrightarrow B^{r_2,s_2}\ot B^{r_1,s_1}
\]
sends $T_1\ot T_2$ to $\td{T}_2\ot \td{T}_1$
if and only if
\[
col(T_2)\rightarrow T_1={col(\td{T}_1)\rightarrow \td{T}_2}.
\]
Moreover, the energy function $H(T_1\ot T_2)$ is given, up to additive constant, 
by the number of nodes in the shape 
of $col(T_2)\rightarrow T_1$ that are strictly to the right of the $\max(s_1,s_2)$-th column.
\end{thm}

\pf
If $T_1\ot T_2$ is a genuine highest weight vector, then the statement is true by 
Lemma \ref{lem:comb R}. Since the insertion algorithm commutes with the action 
of crystal operators $\te_i,\tf_i$ for $i\in \ov{I}$ and each crystal graph $B(\widehat{\la})$ are
connected by $\te_i,\tf_i$ for $i\in \ov{I}$, the former statement follows.

Note that the statement on $H$ immediately follows 
from Lemmas \ref{lem:genuine highest-1} and \ref{lem:genuine highest-2},
\eqref{def H} and the arguments in the proof of Theorem \ref{th:connectedness of tensor product}.
\qed

\begin{rem}{\rm
The description of combinatorial $R$ matrix can be also described in terms of row insertion in a similar way.
}
\end{rem}

\begin{thm}
For $(s_1^{r_1}), (s_2^{r_2}), (s_3^{r_3})\in \cP_{M|N}$, the Yang-Baxter equation
\[
(R\ot1)(1\ot R)(R\ot1)=(1\ot R)(R\ot1)(1\ot R)
\]
holds on $B^{r_1,s_1}\ot B^{r_2,s_2}\ot B^{r_3,s_3}$.
\end{thm}

\pf
As we explained in Section \ref{subsec:comb R} the combinatorial $R$ matrix for
$B^{r_i,s_i}\ot B^{r_j,s_j}$ ($(i,j)=(1,2),(1,3),(2,3)$) is the $q\to0$ limit of the
$R$ matrix for $\W^{(r_i)}_{s_i}\ot\W^{(r_j)}_{s_j}$. Hence the 
Yang-Baxter equation for the $R$ matrix gives rise to that 
for the combinatorial $R$ matrix.
\qed

\begin{rem}
{\rm
Special cases of the combinatorial $R$ matrix or the energy function such as
the $r_1=r_2=1$ case or the $s_1=s_2=1$ case already appeared in literature.
In \cite{KP} the energy function was reinterpreted combinatorially through 
the so-called Howe duality. In \cite{K} a combinatorial rule as in \cite{NY}
was given for the combinatorial $R$ matrix and the energy function by
taking the $q\to0$ limit of the $R$ matrix associated to $\U(\e)$.
}
\end{rem}

\end{document}